\documentclass[11pt,a4paper]{article}

\usepackage{graphicx}
\usepackage{amsfonts}
\usepackage{amsmath}

\makeatletter
\newdimen\mynormalparindent
\def\mymakefnmark{}
\def\mymakefntext{\indent\mymakefnmark}
\long\def\myfootnotetext#1{\insert\footins{%
  \normalfont\footnotesize
  \interlinepenalty\interfootnotelinepenalty
  \splittopskip\footnotesep \splitmaxdepth \dp\strutbox
  \floatingpenalty\@MM \hsize\columnwidth
  \@parboxrestore \parindent\mynormalparindent \sloppy
  \mymakefntext{\rule\z@\footnotesep\ignorespaces#1\unskip\strut\par}}}

\makeatother

\long\def\symbolfootnote[#1]#2{\begingroup%
\def\thefootnote{\fnsymbol{footnote}}\footnote[#1]{#2}\endgroup}

\newtheorem{basic}{Basic}[section]

\newtheorem{lem}[basic]{Lemma}
\newtheorem{propos}[basic]{Proposition}
\newtheorem{thm}[basic]{Theorem}

\newtheorem{cor}[basic]{Corollary}

\newtheorem{con}[basic]{Conjecture}

\newcommand{\bdm}{\begin{displaymath}}
\newcommand{\edm}{\end{displaymath}}
\newcommand{\be}{\begin{equation}}
\newcommand{\ee}{\end{equation}}
\newcommand{\ep}{\vspace{-3mm}\hfill\mbox{$\Box$}\\}

\newcommand{\R}{\mathbb{R}}
\newcommand{\wt}{\widetilde}
\newcommand{\wh}{\widehat}

\def\theequation{\thesection.\arabic{equation}}

\begin{document}

\begin{center}
{\bf\Large On the (non)existence of best low-rank approximations\\ of
generic $I\times J\times 2$ arrays}\\

\vspace{1cm} Alwin Stegeman
\symbolfootnote[2]{A. Stegeman is with the Heijmans Institute for
Psychological Research, University of Groningen, Grote Kruisstraat 2/1,
9712 TS Groningen, The Netherlands, phone: ++31 50 363 6193, fax: ++31 50
363 6304, e-mail: a.w.stegeman@rug.nl, URL:
http://www.gmw.rug.nl/$\sim$stegeman.}\\
\vspace{5mm}
\today
\vspace{5mm}

\begin{abstract}
\noindent Several conjectures and partial proofs have been formulated on
the (non)existence of a best low-rank approximation of
real-valued $I\times J\times 2$ arrays. We analyze this problem using
the Generalized Schur Decomposition and prove (non)existence of a best
rank-$R$ approximation for generic $I\times J\times 2$ arrays, for all
values of $I,J,R$. Moreover, for cases where a best rank-$R$ approximation
exists on a set of positive volume only, we provide easy-to-check necessary
and sufficient conditions for the existence of a best rank-$R$
approximation. \\~\\

\noindent{\em Keywords}: tensor decomposition, low-rank
approximation, Candecomp, Parafac, Generalized Schur Decomposition.\\~\\
\noindent
{\em AMS subject classifications}: 15A18, 15A22, 15A69, 49M27, 62H25.
\end{abstract}\end{center}
\vspace{3cm}
\noindent{\bf Acknowledgments.} This manuscript has benefited from comments of anonymous reviewers who carefully checked the proofs and suggested changes to improve their readability. Also, an anonymous reviewer has stimulated the search for a theoretical proof or rigorous numerical evidence for Conjecture 3.2. Unfortunately, neither have been found by the author. 

This research was supported by the Dutch Organisation for Scientific Research (NWO), VIDI grant 452-08-001.

\newpage
\section{Introduction}
\setcounter{equation}{0}
We consider the problem of finding a best low-rank approximation to a
generic three-way array or order-3 tensor ${\cal Z}\in\R^{I\times J\times
K}$. The rank of a three-way array ${\cal Y}$ is defined
as the smallest number of rank-1 arrays whose sum equals ${\cal Y}$. A
three-way array has rank 1 if it is the outer vector product of three
nonzero vectors. The outer vector product ${\bf Y}={\bf a}\circ{\bf b}={\bf
a}\,{\bf b}^T$ is a rank-1 matrix (or order-2 tensor) with entries
$y_{ij}=a_i\,b_j$. The outer vector product ${\cal Y}={\bf a}\circ{\bf
b}\circ{\bf c}$ is rank-1 tensor with entries $y_{ijk}=a_i\,b_j\,c_k$. The
problem of finding a best rank-$R$ approximation to ${\cal Z}$ can be
denoted as
\be
\label{eq-CPD}
\min_{\substack{{\bf a}_r\in\R^I,\,{\bf b}_r\in\R^J,\,{\bf c}_r\in\R^K, \\
r=1,\ldots,R }}
\quad\|{\cal Z}-\sum_{r=1}^R ({\bf a}_r\circ{\bf b}_r\circ{\bf c}_r)\|^2_F\,,
\ee

\noindent where $\|\cdot\|_F$ denotes the Frobenius norm (i.e., the square
root of the sum-of-squares). For $N$-way arrays (or order-$N$ tensors),
this problem has been introduced by Hitchcock \cite{H1} \cite{H2}. The form
of the rank-$R$ approximation is known as Candecomp/Parafac \cite{Har}
\cite{CC} and also as Canonical Polyadic Decomposition (CPD). It can be
seen as a multi-way (or higher-order) generalization of component analysis
for matrices. Applications of the CPD are found in chemometrics \cite{SBG},
the behavioral sciences \cite{Kro}, signal processing \cite{CoLa}
\cite{DL5}, algebraic complexity theory \cite{Bini-1} \cite{Bini-2} (see
\cite{Ste2} for a discussion), and data mining in general. An overview of
applications of tensor decompositions can be found in \cite{KB}
\cite{Acar}. For the computation of a best low-rank approximation an
iterative algorithm is used. For an overview and comparison of CPD
algorithms, see \cite{Hop} \cite{TB} \cite{CLA}. It has been
proven that determining the rank of an order-3 tensor or computing its
best rank-1 approximation are NP-hard problems \cite{HiLim}.

We denote the frontal $I\times J$ slices of ${\cal Z}\in\R^{I\times
J\times K}$ as ${\bf Z}_k$, $k=1,\ldots,K$. Let ${\bf A}=[{\bf
a}_1\,\ldots\,{\bf a}_R]$, ${\bf B}=[{\bf b}_1\,\ldots\,{\bf b}_R]$, and
${\bf C}=[{\bf c}_1\,\ldots,{\bf c}_R]$. Problem (\ref{eq-CPD})
can be written slicewise as
\be
\label{eq-CPmat}
\min_{\substack{{\bf A}\in\R^{I\times R},\,{\bf B}\in\R^{J\times R},\\
{\rm diag}({\bf C}_k)\in\R^R,\,k=1,\ldots,K}}
\quad \sum_{k=1}^K \|{\bf Z}_k-{\bf
A}\,{\bf C}_k\,{\bf B}^T\|^2_F\,, \ee

\noindent where ${\bf C}_k$ is $R\times R$ diagonal with row $k$ of ${\bf
C}$ as diagonal, $k=1,\ldots,K$. The set of $I\times J\times K$ arrays with
rank at most $R$ is denoted by
\be
S_R(I,J,K)=\{{\cal Y}\in\R^{I\times J\times K}:\;{\rm rank}({\cal Y})\le
R\}\,.
\ee

\noindent Problem (\ref{eq-CPD}) can also be written as:
\be
\label{prob-CPD}
\min_{{\cal Y}\in S_R(I,J,K)} \|{\cal Z}-{\cal Y}\|^2_F\,.
\ee

\noindent Unfortunately, for $R\ge 2$, the problem may not have an optimal
solution because the set $S_R(I,J,K)$ is not closed \cite{DSL}. In such a
case, trying to compute a best rank-$R$ approximation yields a rank-$R$
sequence converging to a boundary point ${\cal X}$ of $S_R(I,J,K)$ with
rank$({\cal X})>R$. As a result, while running the iterative CPD algorithm,
the decrease of the objective function becomes very slow, and some (groups
of) columns of ${\bf A}$, ${\bf B}$, and ${\bf C}$ become nearly linearly
dependent, while their norms increase without bound \cite{KHL} \cite{KDS} \cite{DSL}.
This phenomenon is known as ``diverging CP components" or ``degenerate
solutions" or ``diverging rank-1 terms". Needless to say, diverging rank-1
terms should be avoided if an interpretation of the rank-1 terms is needed.
Note that diverging rank-1 terms are used in algebraic complexity theory to
obtain a fast and arbitrarily accurate approximation to the computation of
bilinear forms (see \cite{Ste2} for a discussion).

Nonexistence of a best rank-$R$ approximation can be avoided by imposing
constraints on the rank-1 terms in $({\bf A},{\bf B},{\bf C})$. Imposing
orthogonality constraints on (one of) the component matrices guarantees
existence of a best rank-$R$ approximation \cite{KDS}, and the same is true
for nonnegative ${\cal Z}$ under the restriction of nonnegative
${\bf A},{\bf B},{\bf C}$ \cite{LiCo}. Also, \cite{LiCo2} show that
constraining the magnitude of the inner products between pairs of columns
of ${\bf A},{\bf B},{\bf C}$ guarantees existence of a best rank-$R$
approximation. However, imposing constraints will not be suitable for all
CPD applications. As an alternative to deal with diverging rank-1 terms,
methods have been developed to obtain the limit point ${\cal X}$ of the
diverging rank-$R$ sequence and a sparse decomposition of ${\cal X}$
\cite{SDL} \cite{RoG} \cite{Ste-limit} \cite{Ste-Jordan} \cite{Ste-TV}.

There are very few theoretical results on the (non)existence of a best
rank-$R$ approximation for specific three-way arrays or sizes. It has been
proven that $2\times 2\times 2$ arrays of rank 3 do not have a best rank-2
approximation \cite{DSL}, and conjectures on $I\times J\times 2$ arrays are
formulated and partly proven in \cite{Ste2}. In simulation studies with
random ${\cal Z}$, diverging rank-1 terms occur very often \cite{Ste}
\cite{Ste2} \cite{Ste3} \cite{Ste-limit}. Although diverging rank-1 terms
may also occur due to a bad choice of starting point for the iterative
algorithm \cite{Paa2} \cite{Ste4}, if trying many random starting points
does not help, then this is strong evidence for nonexistence of a best
rank-$R$ approximation.

In this paper, we consider (non)existence of best rank-$R$ approximations
for generic $I\times J\times 2$ arrays. The use of the term {\em generic}
implies that the entries are randomly sampled from an $IJ2$-dimensional
continuous distribution (for which sets of positive Lebesgue measure also
have positive probability). Properties that hold for a generic array hold
``with probability one", ``almost surely", or ``almost everywhere".
Properties that hold on a set of positive Lebesgue measure but not almost
everywhere, hold on a set of ``positive volume" or ``with positive
probability". Using the relations between the CPD and the Generalized Schur
Decomposition (GSD) formulated in \cite{LMV} \cite{SDL} \cite{Ste-arxiv},
we are able to prove the conjectures formulated in \cite{Ste2}. Our main
result concerns generic $I\times I\times 2$ arrays, which have ranks $I$
and $I+1$ on sets of positive Lebesgue measure. It has been conjectured
that generic $I\times I\times 2$ arrays of rank $I+1$ do not have a best
rank-$I$ approximation \cite{Ste} \cite{Ste2}. So far, this has only been
proven for $I=2$ \cite{DSL}. We provide a proof for $I\ge 2$. Our proofs of
the (non)existence of best rank-$R$ approximations for generic $I\times
J\times 2$ arrays make use of our main result. In some cases, we
prove that existence of a best rank-$R$ approximation holds on a set of
positive volume only. For such arrays we also provide easy-to-check
necessary and sufficient conditions for the existence of a best rank-$R$
approximation.

Classically, $I\times J\times 2$ arrays were classified as matrix pencils,
where a matrix pencil $\mu\,{\bf X}_1+\lambda\,{\bf X}_2$ consists of two
$I\times J$ matrices ${\bf X}_1$ and ${\bf X}_2$ and scalars $\mu$ and
$\lambda$. A matrix pencil is called {\em regular} if both ${\bf X}_1$ and
${\bf X}_2$ are square matrices and there exist $\mu$ and $\lambda$ such
that det$(\mu{\bf X}_1+\lambda{\bf X}_2)\neq 0$. In all other cases, the
pencil is called {\em singular}. For regular matrix pencils, equivalence
results and a canonical form were established by Weierstrass \cite{Wei}.
The corresponding theory for singular pencils was developed by Kronecker
\cite{Kron}. Ja' Ja' \cite{JJ1} extended Kronecker's \cite{Kron}
equivalence results for matrix pencils to $I\times J\times 2$ arrays.
As an anonymous reviewer pointed out, orbit classification results for
matrix pencils such as those recently obtained for complex pencils in
\cite{Perv} could form a different approach to prove
(non)existence of best rank-$R$ approximations. However, we have taken a
different approach via the relation between the CPD and GSD. 
A more detailed discussion of the relation between classical matrix pencil theory
and the rank of real $I\times J\times 2$ arrays can be found in \cite{Ste2}.

The paper is organized as follows. In section 2, we consider the relation
between the CPD and GSD for $I\times J\times 2$ arrays and state the
conjectures of \cite{Ste2}. In section 3, we formulate our main result for
$I\times I\times 2$ arrays and $R=I$, and sketch its proof. The proof
itself is contained in the appendix. In section 4, we prove a case that
cannot be proven by using the GSD. In section 5, we extend our
analysis and proof from section 3 to $I\times J\times 2$ arrays and
$R\le\min(I,J)$. Finally, section 6 contains a discussion of our findings.

We use the following notation. The notation ${\cal Y}$, ${\bf Y}$, ${\bf
y}$, $y$ is used for a three-way array, a matrix, a column vector, and a
scalar, respectively. All arrays, matrices, vectors, and scalars are
real-valued. Matrix transpose and inverse are denoted as ${\bf
Y}^T$ and ${\bf Y}^{-1}$, respectively. An zero matrix of size $p\times
q$ is denoted by ${\bf O}_{p,q}$. An zero column vector is
denoted by ${\bf 0}$. A $p\times p$ matrix ${\bf Y}$ is called orthonormal 
if ${\bf Y}^T{\bf Y}={\bf YY}^T={\bf I}_p$. A $p\times q$ matrix has orthogonal 
columns if ${\bf Y}^T{\bf Y}$ is diagonal.

\section{The CPD and GSD for $I\times J\times 2$ arrays}
\setcounter{equation}{0}
We begin by defining the Generalized Schur Decomposition (GSD) for $I\times
J\times 2$ arrays. Analogous to (\ref{eq-CPmat}), fitting a GSD to ${\cal
Z}$ can be written slicewise as
\be
\label{eq-GSD}
\min_{\substack{{\bf Q}_a\in\R^{I\times R},\,{\bf Q}_b\in\R^{J\times R},\\
{\bf Q}_a^T{\bf Q}_a={\bf Q}_b^T{\bf Q}_b={\bf I}_R,\\
{\bf R}_k\in\R^{R\times R}\;{\rm upper\;triangular},\,k=1,2.}}
\quad \sum_{k=1}^2 \|{\bf Z}_k-{\bf Q}_a\,{\bf R}_k\,{\bf
Q}_b^T\|^2_F\,.
\ee

\noindent Note that the GSD is only defined for $R\le\min(I,J)$.
We define the GSD solution set as
\be
\label{eq-Pset}
P_R(I,J,2)=\{{\cal Y}\in\R^{I\times J\times 2}:\;{\bf Y}_k={\bf Q}_a\,{\bf
R}_k\,{\bf Q}_b^T\,,\;k=1,2\}\,.
\ee

\noindent It has been shown that $P_R(I,J,2)$ is equal to the
closure of $S_R(I,J,2)$ \cite{SDL} \cite{Ste-arxiv}. Moreover, a best
fitting GSD always exists and it can be transformed to a best rank-$R$
approximation if it exists \cite{SDL}. If a best rank-$R$ approximation
does not exist, then a CPD algorithm trying to find a best rank-$R$
approximation yields a sequence of rank-$R$ arrays converging to an optimal
solution of (\ref{eq-GSD}), and the CPD sequence features diverging
components.

Showing that ${\cal Z}$ has no best rank-$R$ approximation is equivalent to
showing that all optimal solutions of (\ref{eq-GSD}) have rank larger than
$R$. Let ${\bf G}_a$ and ${\bf G}_b$ be such that $\wt{\bf Q}_a=[{\bf
Q}_a\;{\bf G}_a]$ and $\wt{\bf Q}_b=[{\bf Q}_b\;{\bf G}_b]$ are square and
orthonormal matrices. When the slices of the GSD solution array are
premultiplied by $\wt{\bf Q}_a^T$ and postmultiplied by $\wt{\bf Q}_b$, we
obtain slices
\be
\left[\begin{array}{cc} {\bf R}_k & {\bf O}_{R,J-R}\\ {\bf O}_{I-R,R} &
{\bf O}_{I-R,J-R}\end{array}\right]\,,\quad\quad k=1,2\,, \ee

\noindent where ${\bf O}_{p,q}$ denotes an zero $p\times q$ matrix.
This implies that the rank of the GSD solution array is equal to the rank
of the $R\times R\times 2$ array ${\cal R}$ with slices ${\bf R}_1$ and
${\bf R}_2$. To establish the rank of ${\cal R}$, we use the following
lemma.

\begin{lem}
\label{lem-rank}
Let ${\cal Y}\in\R^{R\times R\times 2}$ with nonsingular $R\times R$ slices
${\bf Y}_1$ and ${\bf Y}_2$. The following statements hold:
\begin{itemize}
\item[$(i)$] If ${\bf Y}_2{\bf Y}_1^{-1}$ has $R$ real eigenvalues and is
diagonalizable, then ${\cal Y}$ has rank $R$.
\item[$(ii)$] If ${\bf Y}_2{\bf Y}_1^{-1}$ has $R$ real eigenvalues but
is not diagonalizable, then ${\cal Y}$ has at least rank $R+1$.
\item[$(iii)$] If ${\bf Y}_2{\bf Y}_1^{-1}$ has at least one pair of
complex eigenvalues, then ${\cal Y}$ has at least rank $R+1$.
\end{itemize}
\end{lem}

\noindent {\bf Proof.} See \cite[section 3]{JJ}.
\ep

\noindent Suppose that ${\bf R}_1$ and ${\bf R}_2$ are nonsingular. Since
${\bf R}_2{\bf R}_1^{-1}$ is upper triangular, it has $R$ real
eigenvalues. By Lemma~\ref{lem-rank}, the rank of ${\cal R}$ is $R$ when
${\bf R}_2{\bf R}_1^{-1}$ has $R$ linearly independent eigenvectors.
Otherwise, the rank of ${\cal R}$ is larger than $R$. In this case, the GSD
solution array is the limit point of a CPD sequence featuring diverging
rank-1 terms. Moreover, the diverging rank-1 terms are defined by groups of
identical eigenvalues that do not have the same number of linearly
independent associated eigenvectors \cite{Ste} \cite{Ste2} \cite{SDL}
\cite{SDL-TR}.

In this paper, we consider the conjectures of \cite{Ste2} on the
(non)existence of best low-rank approximations for generic $I\times J\times
2$ arrays. These conjectures are given in Table~\ref{tab-1}. Note that
existence of a best rank-$R$ approximation is formulated in terms of
volume, but can analogously be formulated in terms of probability. The rank
values for the generic arrays are derived from the following. For a generic
$I\times I\times 2$ array ${\cal Z}$, the matrix ${\bf Z}_2{\bf Z}_1^{-1}$
has $I$ distinct eigenvalues. By Lemma~\ref{lem-rank} and \cite{TBK}
\cite{Ste}, the array satisfies either $(i)$ and has rank $I$, or $(iii)$
and has rank $I+1$; see also \cite{TBK}. This is also formulated as
$I\times I\times 2$ arrays having {\em typical rank} $\{I,I+1\}$. For
generic $I\times J\times 2$ arrays with $I>J\ge 2$, the rank is given by
$\min(I,2J)$ \cite{TBK}. In other words, $I\times J\times 2$ arrays with
$I>J\ge 2$ have {\em generic rank} $\min(I,2J)$. The notion of typical rank
is used when several rank values occur on sets of positive Lebesgue
measure.

In cases 1, 4, and 6 in Table~\ref{tab-1}, the value of $R$ is larger than
or equal to the rank of ${\cal Z}$. Hence, in these cases the best rank-$R$
approximation of ${\cal Z}$ is ${\cal Z}$ itself. Case 2 is proven in
section 3. Cases 3, 5, 8, and 9 are proven in section 5. In case 7 we
have $R>J$ and cannot use the GSD to analyze the problem. This case is
proven in section 4.

\begin{table}[t]
\begin{center}
\begin{tabular}{|c||c|c|c|c|}
\hline Case & ${\cal Z}\in\R^{I\times J\times 2}$ &
rank$({\cal Z})$ & $R$ & Best rank-$R$ approx. exists ?\\[2mm]

\hline\hline

1 & $I=J$ & $I+1$ & $R\ge I+1$ & always \\[2mm]

2 & $I=J$ & $I+1$ & $R=I$ & zero volume \\[2mm]

3 & $I=J$ & $I+1$ & $R<I$ & positive volume \\[2mm]

\hline 4 & $I=J$ & $I$ & $R\ge I$ & always \\[2mm]

5 & $I=J$ & $I$ & $R<I$ & positive volume \\[2mm]

\hline 6 & $I>J$ & $\min(I,2J)$ & $R\ge\min(I,2J)$ & always \\[2mm]

7 & $I>J$ & $\min(I,2J)$ & $\min(I,2J)>R>J$ & almost everywhere \\[2mm]

8 & $I>J$ & $\min(I,2J)$ & $R=J$ & positive volume \\[2mm]

9 & $I>J$ & $\min(I,2J)$ & $R<J$ & positive volume \\[2mm]

\hline

\end{tabular}
\end{center}

\caption{Results (cases 1, 4, and 6) and conjectures (cases 2, 3, 5, 7, 8,
and 9) of \cite{Ste2} on the existence of a best rank-$R$ approximation of
generic $I\times J\times 2$ arrays. Here, $I\ge J\ge 2$ and $R\ge 2$.}
\label{tab-1} \end{table}\vspace{5mm}

\section{Case 2: $I\times I\times 2$ arrays of rank $I+1$ and $R=I$}
\setcounter{equation}{0}
We consider the GSD problem for generic $I\times I\times 2$ arrays
and $R=I$. We rewrite the GSD problem (\ref{eq-GSD}) as
\be
\label{eq-GSD2}
\min_{\substack{{\bf Q}_a\in\R^{I\times R},\,{\bf Q}_b\in\R^{J\times R},\\
{\bf Q}_a^T{\bf Q}_a={\bf Q}_b^T{\bf Q}_b={\bf I}_R,\\
{\bf R}_k\in\R^{R\times R}\;{\rm upper\;triangular},\,k=1,2.}}
\quad \sum_{k=1}^2 \|{\bf Q}_a^T\,{\bf Z}_k\,{\bf Q}_b-{\bf
R}_k\|^2_F\,. \ee

\noindent For each ${\bf Q}_a$ and ${\bf Q}_b$, the optimal ${\bf R}_k$
are found as the upper triangular parts of ${\bf Q}_a^T{\bf Z}_k{\bf Q}_b$,
respectively. Hence, problem (\ref{eq-GSD2}) can be written as
\be
\label{eq-GSD3}
\min_{\substack{{\bf Q}_a\in\R^{I\times R},\,{\bf Q}_b\in\R^{J\times R},\\
{\bf Q}_a^T{\bf Q}_a={\bf Q}_b^T{\bf Q}_b={\bf I}_R}}
\quad \sum_{k=1}^2 \|{\bf Q}_a^T\,{\bf Z}_k\,{\bf
Q}_b\|^2_{LFs}\,,
\ee

\noindent where $\|\cdot\|_{LFs}$ denotes the Frobenius norm of the strictly
lower triangular part. The optimal ${\bf Q}_a$ and ${\bf Q}_b$ can be
obtained by iterating over Givens rotations (De Lathauwer, De Moor, and
Vandewalle \cite{LMV}). The optimal ${\bf Q}_a$ and ${\bf Q}_b$ are then
the products of the consecutive optimal Givens rotation matrices. Each
rotation affects rows and columns $i$ and $j$ ($i<j$) of ${\bf Q}_a^T{\bf
Z}_k{\bf Q}_b$, $k=1,2$. For rotation $(i,j)$, the corresponding Givens
rotation matrices ${\bf U}_a$ and ${\bf U}_b$ are equal to ${\bf I}_I$
except:
\be
\label{eq-rot-Ua}
({\bf U}_a)_{ii}=({\bf U}_a)_{jj}=\cos(\alpha)\,,\quad\quad\quad
({\bf U}_a)_{ji}=-({\bf U}_a)_{ij}=\sin(\alpha)\,,
\ee
\be
\label{eq-rot-Ub}
({\bf U}_b)_{ii}=({\bf U}_b)_{jj}=\cos(\beta)\,,\quad\quad\quad
({\bf U}_b)_{ji}=-({\bf U}_b)_{ij}=\sin(\beta)\,.
\ee

\noindent The Jacobi-type algorithm of \cite{LMV} to solve problem
(\ref{eq-GSD3}) iterates over all rotations $(i,j)$, $1\le i<j\le I$. In
each iteration, $\alpha$ and $\beta$ are computed that minimize
$\sum_{k=1}^2 \|{\bf U}_a^T\,{\bf Q}_a^T\,{\bf Z}_k\,{\bf Q}_b\,{\bf
U}_b\|^2_{LFs}$, where ${\bf Q}_a$ and ${\bf Q}_b$ are the current updates.
Next, ${\bf Q}_a$ is replaced by ${\bf Q}_a{\bf U}_a$ and ${\bf Q}_b$ is
replaced by ${\bf Q}_b{\bf U}_b$. A necessary condition for reaching an
optimal solution is that no rotation $(i,j)$ can further decrease the
objective function in (\ref{eq-GSD3}). To derive the equations defining
local minima for each rotation $(i,j)$, we use the following lemma.

\begin{lem}
\label{lem-statorth}
For vectors ${\bf x},{\bf y}\in\R^p$ and $\alpha\in\R$, define the rotation
\bdm
[\tilde{\bf x}\;\;\tilde{\bf y}]=[{\bf x}\;\;{\bf
y}]\;\left[\begin{array}{cc}
\cos(\alpha) & -\sin(\alpha) \\
\sin(\alpha) & \cos(\alpha)\end{array}\right]\,.
\edm

\noindent Let $f(\alpha)=\|\tilde{\bf x}\|^2$. We have
\be
\label{eq-derivatives}
\frac{\partial f}{\partial\alpha}=2\,\tilde{\bf x}^T\tilde{\bf y}\,,
\quad\quad\quad
\frac{\partial^2 f}{\partial\alpha^2}=
2\,(\tilde{\bf y}^T\tilde{\bf y}-\tilde{\bf x}^T\tilde{\bf x})\,.
\ee

\noindent Moreover, if $\frac{\partial f}{\partial\alpha}
(\alpha)=\frac{\partial^2 f}{\partial\alpha^2}(\alpha)=0$ for some
$\alpha$, then $f(\alpha)= {\bf x}^T{\bf x}$ is constant.
\end{lem}

\noindent {\bf Proof.} We write $\|\tilde{\bf x}\|^2=\sum_{i=1}^p
(\cos(\alpha)\,x_i+\sin(\alpha)\,y_i)^2$. The first derivative
in (\ref{eq-derivatives}) follows from
\bdm
\frac{\partial f}{\partial\alpha}=
2\,\sum_{i=1}^p
(\cos(\alpha)\,x_i+\sin(\alpha)\,y_i)\,(-\sin(\alpha)\,x_i+\cos(\alpha)\,y_i
)= 2\,\sum_{i=1}^p \tilde{x}_i\,\tilde{y}_i = 2\,\tilde{\bf
x}^T\tilde{\bf y}\,. \edm

\noindent The second derivative is obtained as
\begin{eqnarray*}
\frac{\partial^2 f}{\partial\alpha^2} &=& 2\,\frac{\partial f}{\partial\alpha}
\sum_{i=1}^p
(\cos(\alpha)\,x_i+\sin(\alpha)\,y_i)\,(-\sin(\alpha)\,x_i+\cos(\alpha)\,y_i)
\\
&=& 2\,\sum_{i=1}^p ((-\sin(\alpha)\,x_i+\cos(\alpha)\,y_i)^2 -
(\cos(\alpha)\,x_i+\sin(\alpha)\,y_i)^2) \\
&=& 2\,\sum_{i=1}^p \tilde{y}_i^2-\tilde{x}_i^2=
2\,(\tilde{\bf y}^T\tilde{\bf y}-\tilde{\bf x}^T\tilde{\bf x})\,.
\end{eqnarray*}

\noindent Next, suppose the first and second derivatives are zero
for some $\alpha$. That is, $\tilde{\bf x}^T\tilde{\bf y}=0$ and
$\tilde{\bf x}^T\tilde{\bf x}=\tilde{\bf y}^T\tilde{\bf y}$ for some
$\alpha$. We write
\begin{eqnarray}
\label{eq-first-zero}
\tilde{\bf x}^T\tilde{\bf y} &=& (\cos^2(\alpha)-\sin^2(\alpha))\,
{\bf x}^T{\bf y}+\sin(\alpha)\cos(\alpha)\,({\bf y}^T{\bf y}-
{\bf x}^T{\bf x})=0\,, \\
\label{eq-second-zero}
\tilde{\bf y}^T\tilde{\bf y}-\tilde{\bf x}^T\tilde{\bf x} &=&
(\cos^2(\alpha)-\sin^2(\alpha))\,({\bf y}^T{\bf y}-{\bf x}^T{\bf x})-
4\,\sin(\alpha)\cos(\alpha)\,{\bf x}^T{\bf y}=0\,, \\
\label{eq-f}
f(\alpha)=\tilde{\bf x}^T\tilde{\bf x} &=& \cos^2(\alpha)\,
{\bf x}^T{\bf x}+\sin^2(\alpha)\,{\bf y}^T{\bf y}+2\,\sin(\alpha)
\cos(\alpha)\,{\bf x}^T{\bf y}\,.
\end{eqnarray}

\noindent When $\sin(\alpha)=0$ or $\cos(\alpha)=0$, it follows from
(\ref{eq-first-zero})-(\ref{eq-second-zero}) that ${\bf x}^T{\bf y}=0$ and
${\bf x}^T{\bf x}={\bf y}^T{\bf y}$. By (\ref{eq-f}), this implies
the desired result $f(\alpha)={\bf x}^T{\bf x}$. Next, suppose
$\sin(\alpha)\cos(\alpha)\neq 0$. Combining
(\ref{eq-first-zero})-(\ref{eq-second-zero}) yields
\be
\label{eq-comb}
{\bf y}^T{\bf y}-{\bf x}^T{\bf x}=-\left(\frac{(\cos^2(\alpha)
-\sin^2(\alpha))^2}{4\,\sin^2(\alpha)\cos^2(\alpha)}\right)\;
({\bf y}^T{\bf y} -{\bf x}^T{\bf x})\,.
\ee

\noindent Since the term depending on $\alpha$ in (\ref{eq-comb}) is
nonpositive, it follows that ${\bf x}^T{\bf x}={\bf y}^T{\bf y}$. Then
${\bf x}^T{\bf y}=0$ follows from
(\ref{eq-first-zero})-(\ref{eq-second-zero}), and we again obtain
$f(\alpha)={\bf x}^T{\bf x}$. This completes the proof.\ep

\noindent Let $\wt{\bf Z}_k={\bf Q}_a^T\,{\bf Z}_k\,{\bf Q}_b$, $k=1,2$, and
define the 2-dimensional vectors
\be
\label{eq-tildez}
\tilde{\bf z}_{(m,n)}=\left(\begin{array}{c}
(\wt{\bf Z}_1)_{mn}\\
(\wt{\bf Z}_2)_{mn}\end{array}\right)=\left(\begin{array}{c}
\tilde{z}_{mn1}\\ \tilde{z}_{mn2}\end{array}\right)\,,\quad\quad\quad
m=1,\ldots,I\,,\quad n=1,\ldots,I\,.
\ee

\noindent The vectors (\ref{eq-tildez}) are also known as mode-3 vectors of the array 
with slices $\wt{\bf Z}_1$ and $\wt{\bf Z}_2$. As in \cite{LMV}, we determine the stationary points for
rotation $(i,j)$ by setting the derivatives with respect to $\alpha$ and
$\beta$ of $\sum_{k=1}^2 \|{\bf U}_a^T\,\wt{\bf Z}_k\,{\bf U}_b\|^2_{LFs}$
equal to zero. When rotating rows $i$ and $j$ (with $i<j$) the
entries $(i,r)$ and $(j,r)$ with $r=1,\ldots,i-1$ stay in the strictly lower
triangular part. Their Frobenius norm is not changed. Analogously, the
entries $(i,r)$ and $(j,r)$ with $r=j,\ldots,I$ stay in the upper
triangular part, and do not affect the objective function (\ref{eq-GSD3}).
Hence, the rotation of rows $i$ and $j$ can change the objective function
only via entries $(i,r)$ and $(j,r)$ with $r=i,\ldots,j-1$. Let $\tilde{\bf y}_k=
[\tilde{z}_{i,i,k},\ldots,\tilde{z}_{i,j-1,k}]^T$ and $\tilde{\bf x}_k=
[\tilde{z}_{j,i,k},\ldots,\tilde{z}_{j,j-1,k}]^T$, $k=1,2$. Next we apply Lemma~\ref{lem-statorth} 
with objective function $f(\alpha)=\|\tilde{\bf x}_1\|^2+\|\tilde{\bf x}_2\|^2$, which yields 
the first-order condition $\tilde{\bf y}_1^T\tilde{\bf x}_1+\tilde{\bf y}_2^T
\tilde{\bf x}_2=0$. We rewrite this condition for a stationary point as
\be
\label{eq-stat1}
\sum_{r=i}^{j-1} \tilde{\bf z}_{(i,r)}^T\tilde{\bf z}_{(j,r)}=0\,,
\quad\quad 1\le i<j\le I\,.
\ee

\noindent When rotating columns $i$ and $j$ (with $i<j$) we obtain
analogously that the objective function can be changed only via entries
$(r,i)$ and $(r,j)$ with $r=i+1,\ldots,j$. Analogous to (\ref{eq-stat1}), 
it follows from Lemma~\ref{lem-statorth} that a stationary point satisfies
\be
\label{eq-stat2}
\sum_{r=i+1}^{j} \tilde{\bf z}_{(r,i)}^T\tilde{\bf z}_{(r,j)}=0\,,
\quad\quad 1\le i<j\le I\,.
\ee

\noindent Equations (\ref{eq-stat1}) and (\ref{eq-stat2}) are the first-order optimality conditions. Hence, in an optimal solution of problem (\ref{eq-GSD3}) equations (\ref{eq-stat1}) and (\ref{eq-stat2}) will hold.

We obtain second-order optimality conditions from Lemma~\ref{lem-statorth},
where positive second derivatives are required for a local minimum for each
rotation. Lemma~\ref{lem-statorth} shows that a second derivative being
zero at a stationary point implies a constant objective function and
infinitely many optimal rotation angles. For the rotation of rows $i$ and $j$ with
$i<j$, we define $\tilde{\bf x}_k$ and $\tilde{\bf y}_k$ as above, $k=1,2$.
For the minimization of $f(\alpha)=\|\tilde{\bf x}_1\|^2+\|\tilde{\bf x}_2\|^2$, the second-order condition 
is given by $\|\tilde{\bf y}_1\|^2+\|\tilde{\bf y}_2\|^2\ge
\|\tilde{\bf x}_1\|^2+\|\tilde{\bf x}_2\|^2$; see Lemma~\ref{lem-statorth}. 
We rewrite this second-order condition as
\be
\label{eq-second1}
\sum_{r=i}^{j-1}\tilde{\bf z}_{(i,r)}^T\tilde{\bf z}_{(i,r)}-
\sum_{r=i}^{j-1}\tilde{\bf z}_{(j,r)}^T\tilde{\bf z}_{(j,r)}\ge 0\,,
\quad\quad 1\le i<j\le I\,.
\ee

\noindent Analogously, the second-order condition for the rotation of columns $i$ and
$j$ with $i<j$ equals
\be
\label{eq-second2}
\sum_{r=i+1}^{j}\tilde{\bf z}_{(r,j)}^T\tilde{\bf z}_{(r,j)}-
\sum_{r=i+1}^{j}\tilde{\bf z}_{(r,i)}^T\tilde{\bf z}_{(r,i)}\ge 0\,,
\quad\quad 1\le i<j\le I\,.
\ee

\noindent We work under the following conjecture of positive second-order conditions.
Numerical evidence and more theoretical context are provided in Appendix A.

\begin{con}
\label{conjec}
Let ${\cal Z}\in\R^{I\times I\times 2}$ be generic with {\rm rank}$({\cal
Z})=I+1$. Let $({\bf Q}_a,{\bf Q}_b,{\bf R}_1,{\bf R}_2)$ be an optimal
solution of the GSD problem $(\ref{eq-GSD})$ with $R=I$. Then the second-order conditions $(\ref{eq-second1})$-$(\ref{eq-second2})$ are strictly positive.
\ep
\end{con}

\noindent We use the first and second-order optimality conditions to obtain
the following result.

\begin{thm}
\label{t-1}
Let ${\cal Z}\in\R^{I\times I\times 2}$ be generic with {\rm rank}$({\cal
Z})=I+1$. Let $({\bf Q}_a,{\bf Q}_b,{\bf R}_1,{\bf R}_2)$ be an optimal
solution of the GSD problem $(\ref{eq-GSD})$ with $R=I$ and strictly positive second-order conditions $(\ref{eq-second1})$-$(\ref{eq-second2})$. Then the rank of
the $I\times I\times 2$ array ${\cal R}$ with slices ${\bf R}_1$ and ${\bf
R}_2$ is larger than $I$.
\end{thm}

\noindent {\bf Proof.} See Appendix B.\ep

\noindent Theorem~\ref{t-1} and Conjecture~\ref{conjec} imply that any optimal solution array of the GSD problem (\ref{eq-GSD}), which has slices ${\bf Q}_a{\bf R}_k{\bf
Q}_b^T$, $k=1,2$, has rank larger than $I$. As explained in section 2, this
is equivalent to ${\cal Z}$ not having a best rank-$I$ approximation.
Hence, we obtain the following.

\begin{cor}
\label{c-1}
Let ${\cal Z}\in\R^{I\times I\times 2}$ be generic with {\rm rank}$({\cal
Z})=I+1$. Then ${\cal Z}$ does not have a best rank-$I$ approximation.\ep
\end{cor}

\noindent Note that the formulation of Corollary~\ref{c-1} is equivalent
to an $I\times I\times 2$ array of rank $I+1$ having a best
rank-$I$ approximation at most on a set of zero volume, as is stated in case 2 of
Table~\ref{tab-1}.

\section{Case 7: $I\times J\times 2$ arrays with $I>J$ and $\min(I,2J)>R>J$}
\setcounter{equation}{0}
Since $R>J$, we cannot use the GSD in this case. We define the set
\be
\label{eq-W}
W_R(I,J,2)=\{{\cal Y}\in\R^{I\times J\times 2}:\;{\rm rank}([{\bf
Y}_1\,|\,{\bf Y}_2])\le R\}\,.
\ee

\noindent and the problem
\be
\label{prob-W}
\min_{{\cal Y}\in W_R(I,J,2)} \|{\cal Z}-{\cal Y}\|^2_F\,. 
\ee

\noindent Stegeman \cite{Ste2} shows that $W_R(I,J,2)$ is the closure of
the rank-$R$ set $S_R(I,J,2)$ when $I>J\ge 2$ and $\min(I,2J)>R>J$. In our
proof below, the set $W_R(I,J,2)$ plays the role of the GSD solution set
$P_R(I,J,2)$ in section 3. We have the following result.

\begin{thm}
\label{t-5}
Let ${\cal Z}\in\R^{I\times J\times 2}$ with $I>J\ge 2$ and
$\min(I,2J)>R>J$ be generic. Then the optimal solution ${\cal X}$ of
problem $(\ref{prob-W})$ is unique and {\rm rank}$({\cal X})=R$.
\end{thm}

\noindent {\bf Proof.} Problem (\ref{prob-W}) is in fact a matrix problem.
Namely, the closest rank-$R$ matrix ${\bf Y}=[{\bf Y}_1\,|\,{\bf Y}_2]$ to
a generic $I\times 2J$ matrix ${\bf Z}=[{\bf Z}_1\,|\,{\bf Z}_2]$ is asked
for. It is well known that this problem is solved by the truncated singular
value decomposition (SVD) of ${\bf Z}$ \cite{EY}. Let the SVD of ${\bf Z}$
be given as ${\bf Z}={\bf U}\,{\bf S}\,{\bf V}^T$. Without loss of
generality we assume $I\ge 2J$. Matrix ${\bf U}$ is $I\times 2J$ and
columnwise orthonormal, ${\bf S}$ is $2J\times 2J$ diagonal and
nonsingular, and ${\bf V}$ is $2J\times 2J$ and orthonormal. The singular
values on the diagonal of ${\bf S}$ are assumed to be in decreasing order.
Since the singular values of ${\bf Z}$ are distinct (${\bf Z}$ is generic),
matrix problem (\ref{prob-W}) has a unique solution ${\bf X}={\bf U}_R\,{\bf
S}_R\,{\bf V}_R^T$. Here, ${\bf U}_R$ and ${\bf V}_R$ contain the first $R$
columns of ${\bf U}$ and ${\bf V}$, respectively, and ${\bf S}_R$ is
$R\times R$ diagonal and contains the $R$ largest singular values.

The optimal solution ${\cal X}$ of problem (\ref{prob-W}) has slices ${\bf
X}_1={\bf U}_R\,{\bf S}_R\,{\bf V}_{R,1}^T$ and ${\bf X}_2={\bf U}_R\,{\bf
S}_R\,{\bf V}_{R,2}^T$, where ${\bf V}_{R,1}^T$ contains columns
$1,\ldots,J$ of ${\bf V}_R^T$, and ${\bf V}_{R,2}^T$ contains columns
$J+1,\ldots,2J$ of ${\bf V}_R^T$. The rank of ${\cal X}$ is equal to the
rank of the $R\times J\times 2$ array ${\cal V}_R$ with slices ${\bf
S}_R{\bf V}_{R,1}^T$ and ${\bf S}_R{\bf V}_{R,2}^T$. Hence, the proof is
complete if we show that rank$({\cal V}_R)=R$.

We have the eigendecomposition ${\bf Z}^T{\bf Z}={\bf V}\,{\bf S}^2\,{\bf
V}^T$, where ${\bf Z}^T{\bf Z}$ is a generic symmetric positive definite $2J\times 2J$ matrix. Since the set of all ${\bf Z}^T{\bf Z}$ has positive measure in the set of symmetric $2J\times 2J$ matrices, it has dimensionality $2J(2J+1)/2$. The set with the parameterization ${\bf V}\,{\bf S}\,{\bf V}^T$ must have the same dimensionality. The dimensionality of the set of all ${\bf S}$ equals $2J$ and the dimensionality of the set of all $2J\times 2J$ orthonormal ${\bf V}$ equals $2J(2J-1)/2$, with the sum being $2J(2J+1)/2$. Hence, ${\bf V}$ may be considered a generic $2J\times 2J$ orthonormal matrix. Analogously, ${\bf V}_R$ may be considered a generic $2J\times R$
columnwise orthonormal matrix. The $R\times J\times 2$ array ${\cal V}_R$
may be considered generic under the condition that the rows of its matrix
unfolding ${\bf S}_R{\bf V}_R^T=[{\bf S}_R{\bf V}_{R,1}^T\,|\,{\bf S}_R{\bf
V}_{R,2}^T]$ are orthogonal. Premultiplying the slices of ${\cal V}_R$ by
a generic $R\times R$ matrix yields a generic $R\times J\times 2$ array,
with rank equal to rank$({\cal V}_R)$. Hence, rank$({\cal V}_R)$ is equal
to the rank of generic $R\times J\times 2$ arrays. When $2J>R>J\ge 2$, the
latter rank is given by $\min(R,2J)=R$ \cite{TBK}. This completes the
proof. \ep

\noindent Since $W_R(I,J,2)$ is the closure of $S_R(I,J,2)$, it
follows that the optimal solution ${\cal X}$ in Theorem~\ref{t-5} is an
optimal solution of the best rank-$R$ approximation problem
(\ref{prob-CPD}). Hence, we obtain the following.

\begin{cor}
\label{c-5}
Let ${\cal Z}\in\R^{I\times J\times 2}$ with $I>J\ge 2$ and
$\min(I,2J)>R>J$ be generic. Then ${\cal Z}$ has a best rank-$R$
approximation.
\ep
\end{cor}

\noindent  Note that the formulation of Corollary~\ref{c-5} is equivalent
to an $I\times J\times 2$ array with $I>J\ge 2$ and $\min(I,2J)>R>J$ having
a best rank-$I$ approximation almost everywhere, as is stated in case
7 of Table~\ref{tab-1}.

\section{Extension to $I\times J\times 2$ arrays and $R\le\min(I,J)$}
\setcounter{equation}{0}
Here, we consider the GSD problem (\ref{eq-GSD}) for cases 3, 5, 8, and 9
in Table~\ref{tab-1}. Hence, we have $R<I$ or $R<J$ or both. Also, $R<$
rank$({\cal Z})$. In these cases, \cite{Ste2} conjectures that the set
of arrays that have a best rank-$R$ approximation, and the set of arrays
that do not have a best rank-$R$ approximation, both have positive volume. 
Below, we analyze this using the GSD framework. In section 5.1, we
consider the GSD algorithm when $R<I$ or $R<J$ or both, which was presented
in \cite{SDL}. We derive equations defining a stationary point, which we
use in our proofs. In section 5.2 we prove case 8, in which $R=J<I$. In
section 5.3 we prove cases 3, 5, and 9, in which $R<\min(I,J)$.

\subsection{The GSD algorithm when $R<I$ or $R<J$ or both}
For $R=I=J$, the optimal ${\bf Q}_a$ and ${\bf Q}_b$ are found by
minimizing the Frobenius norm of the strictly lower triangular parts of
${\bf Q}_a^T\,{\bf Z}_k\,{\bf Q}_b$, $k=1,2$; see (\ref{eq-GSD3}). Since
${\bf Q}_a$ and ${\bf Q}_b$ are orthonormal, we have $\|{\bf Q}_a^T\,{\bf
Z}_k\,{\bf Q}_b\|_F=\|{\bf Z}_k\|_F$, $k=1,2$. This implies that solving
(\ref{eq-GSD3}) is equivalent to maximizing the Frobenius norm of the upper
triangular parts of ${\bf Q}_a^T\,{\bf Z}_k\,{\bf Q}_b$, $k=1,2$.
Analogously, for $R<I$ or $R<J$ or both, we maximize the upper
triangular part of the first $R$ rows and columns of $\wt{\bf Q}_a^T\,{\bf
Z}_k\,\wt{\bf Q}_b$, $k=1,2$. Here, $\wt{\bf Q}_a$ $(I\times I)$ and
$\wt{\bf Q}_b$ $(J\times J)$ are orthonormal, and ${\bf Q}_a$ and ${\bf
Q}_b$ are taken as the first $R$ columns from $\wt{\bf Q}_a$ and $\wt{\bf
Q}_b$, respectively.

Updating $\wt{\bf Q}_a$ and $\wt{\bf Q}_b$ is done via Givens rotations, as
for $R=I=J$ in section 3. We have four different kinds of Givens rotations.
Rotations of rows $i$ and $j$ or columns $i$ and $j$ with $1\le i<j\le R$
are the same as described in section 3. Conditions for stationary points
with respect to these rotations are given by (\ref{eq-stat1}) and
(\ref{eq-stat2}) for $1\le i<j\le R$, where we now define $\wt{\bf
Z}_k=\wt{\bf Q}_a^T\,{\bf Z}_k\,\wt{\bf Q}_b$, $k=1,2$. For convenience, we
repeat these equations as
\be
\label{eq-stat12}
\sum_{r=i}^{j-1}
\tilde{\bf z}_{(i,r)}^T\tilde{\bf z}_{(j,r)}=0\,,\quad\quad
1\le i<j\le R\,,
\ee
\be
\label{eq-stat22}
\sum_{r=i+1}^{j}
\tilde{\bf z}_{(r,i)}^T\tilde{\bf z}_{(r,j)}=0\,,\quad\quad
1\le i<j\le R\,.
\ee

\noindent When $R<I$, we have additional rotations of rows $i$ and $j$ with
$i>R$ or $j>R$ or both. Rotations of rows $i$ and $j$ with $R<i<j$ do not
change the upper triangular part of the first $R$ rows. Hence, they can be
left out of consideration. Rotations of rows $i$ and $j$ with $1\le i\le R$
and $R+1\le j\le I$ change the upper triangular part of the first $R$ rows
via entries $(i,r)$ with $r=i,\ldots,R$. Analogous to (\ref{eq-stat12}) and
(\ref{eq-stat22}), this yields the following equations for stationary
points:
\be
\label{eq-stat3}
\sum_{r=i}^R\tilde{\bf z}_{(i,r)}^T\tilde{\bf z}_{(j,r)}=0\,,\quad\quad
1\le i\le R\,,\quad R+1\le j\le I\,.
\ee

\noindent When $R<J$, we also have rotations of columns $i$ and $j$ with
$i>R$ or $j>R$ or both. Analogous to row rotations, we only need to
consider $i$ and $j$ with $1\le i\le R$ and $R+1\le j\le J$. In the upper
triangular part of the first $R$ columns only the entries $(r,i)$ with
$r=1,\ldots,i$ are changed. This yields the following equations for
stationary points:
\be
\label{eq-stat4}
\sum_{r=1}^i\tilde{\bf z}_{(r,i)}^T\tilde{\bf z}_{(r,j)}=0\,,\quad\quad
1\le i\le R\,,\quad R+1\le j\le J\,.
\ee

\noindent Hence, stationary points of the GSD problem (\ref{eq-GSD})
satisfy (\ref{eq-stat12})--(\ref{eq-stat4}).

For fixed $i\in\{1,\ldots,R\}$, the GSD algorithm of \cite{SDL} combines
row rotations $(i,j)$ for all $j=R+1,\ldots,I$ using a singular value
decomposition. The same holds for column rotations $(i,j)$ with fixed
$i\in\{1,\ldots,R\}$ and all $j=R+1,\ldots,J$. However, the GSD algorithm
can also be programmed in the way described above, i.e., solving each
rotation separately. For each rotation, the optimal rotation angle $\alpha$
can be computed by setting the derivative in (\ref{eq-derivatives}) equal
to zero. After dividing by $\cos^2(\alpha)$, this yields a second degree
polynomial in $\tan(\alpha)$. Numerical experiments show that, for the
same generic ${\cal Z}$, the two GSD algorithms yield different ${\bf
Q}_a$, ${\bf Q}_b$, ${\bf R}_1$, and ${\bf R}_2$, but the GSD solution
array is identical, and also the eigenvalues and number of eigenvectors of
${\bf R}_2{\bf R}_1^{-1}$ are identical.

\subsection{Case 8: $I\times J\times 2$ arrays with $I>J=R$}
We proceed analogous to case 7 in section 4. We define the set $W_R(I,J,2)$
as in (\ref{eq-W}) and consider the best approximation of ${\cal Z}$ from
$W_R(I,J,2)$ in (\ref{prob-W}). We have $S_R(I,J,2)\subset W_R(I,J,2)$, see
\cite{Ste2}. Hence, if the best approximation ${\cal X}$ from the set
$W_R(I,J,2)$ has rank at most $R$, then ${\cal Z}$ has a best rank-$R$
approximation. As in the proof of Theorem~\ref{t-5}, the best approximation
from $W_R(I,J,2)$ is unique and given by the truncated singular value
decomposition (SVD) of ${\bf Z}=[{\bf Z}_1\,|\,{\bf Z}_2]$, which we denote
as ${\bf X}={\bf U}_R\,{\bf S}_R\,{\bf V}_R^T$. The corresponding array
${\cal X}$ has slices ${\bf X}_1={\bf U}_R\,{\bf S}_R\,{\bf V}_{R,1}^T$ and
${\bf X}_2={\bf U}_R\,{\bf S}_R\,{\bf V}_{R,2}^T$, where ${\bf V}_R^T=[
{\bf V}_{R,1}^T\,|\,{\bf V}_{R,2}^T]$. The rank of ${\cal X}$ is equal to
the rank of the $R\times R\times 2$ array ${\cal V}_R$ with $R\times R$
slices ${\bf S}_R{\bf V}_{R,1}^T$ and ${\bf S}_R{\bf V}_{R,2}^T$. As stated
above, rank$({\cal X})=$ rank$({\cal V}_R)\le R$ implies that ${\cal Z}$ has
a best rank-$R$ approximation. The rank of ${\cal V}_R$ can be checked by
making use of Lemma~\ref{lem-rank}.

We have the following result for the case where rank$({\cal X})=$
rank$({\cal V}_R)>R$.

\begin{thm}
\label{t-8}
Let ${\cal Z}\in\R^{I\times J\times 2}$ with $I>J=R\ge 2$ be generic. Let
${\cal X}$ be the optimal solution of problem $(\ref{prob-W})$. Let $({\bf
Q}_a,{\bf Q}_b,{\bf R}_1,{\bf R}_2)$ be an optimal solution of the GSD
problem $(\ref{eq-GSD})$ with strictly positive second-order conditions $(\ref{eq-second1})$-$(\ref{eq-second2})$ for $1\le i<j\le R$. If {\rm rank}$({\cal X})>R$, then the rank of the $R\times R\times 2$ array ${\cal R}$ with slices ${\bf R}_1$ and ${\bf
R}_2$ is larger than $R$.
\end{thm}

\noindent {\bf Proof.} See Appendix C.\ep

\noindent Analogous to Conjecture~\ref{conjec}, we assume that the second-order conditions (\ref{eq-second1})-(\ref{eq-second2}), $1\le i<j\le R$, are strictly positive for generic ${\cal Z}\in\R^{I\times J\times 2}$ with $I>J=R\ge 2$. As in Appendix A, we have found no counterexamples in numerical experiments (results not reported). Analogous to Corollary~\ref{c-1} following from Theorem~\ref{t-1}, we obtain the following.

\begin{cor}
\label{c-8}
Let ${\cal Z}\in\R^{I\times J\times 2}$ with $I>J=R\ge 2$ be generic, and
let ${\cal X}$ be the optimal solution of problem $(\ref{prob-W})$.
If {\rm rank}$({\cal X})>R$, then ${\cal Z}$ does not have a best rank-$R$
approximation.
\ep
\end{cor}

\noindent Corollary~\ref{c-8} implies that we now have an
easy-to-check criterion to determine whether ${\cal Z}$ has a best rank-$R$
approximation or not. First, compute the truncated SVD of ${\bf Z}=[{\bf
Z}_1\,|\,{\bf Z}_2]$ as ${\bf X}={\bf U}_R\,{\bf S}_R\,{\bf V}_R^T$. As in
the proof of Theorem~\ref{t-5}, the array ${\cal V}_R$ with slices ${\bf
S}_R{\bf V}_{R,k}^T$, $k=1,2$, may be considered a generic $R\times R\times
2$ array. Hence, its rank is either $R$ or $R+1$, both on sets of positive
Lebesgue measure \cite{TBK}. Next, compute the eigenvalues of ${\bf
S}_R{\bf V}_{R,2}^T({\bf S}_R{\bf V}_{R,1}^T)^{-1}$ (or just ${\bf
V}_{R,2}^T({\bf V}_{R,1}^T)^{-1}$), which are distinct. If all eigenvalues
are real, then rank$({\cal X})=$ rank$({\cal V}_R)=R$
(Lemma~\ref{lem-rank}) and ${\cal Z}$ has a best rank-$R$ approximation,
which can be taken equal to ${\cal X}$. If some eigenvalues are complex,
then rank$({\cal X})=$ rank$({\cal V}_R)=R+1$ (Lemma~\ref{lem-rank}) and
${\cal Z}$ does not have a best rank-$R$ approximation. Since both
situations occur on sets of positive Lebesgue measure, this completes the
proof of case 8 of Table~\ref{tab-1}.

\subsection{Cases 3, 5, 9: $I\times J\times 2$ arrays with $R<\min(I,J)$}
We proceed analogous to case 8 in section 5.2, except now the situation is
more complicated. We define the set
\be
\label{eq-W359}
\wt{W}_R(I,J,2)=\{{\cal Y}\in\R^{I\times J\times 2}:\;{\rm rank}([{\bf
Y}_1\,|\,{\bf Y}_2])\le R\,,\;{\rm and\;rank}\left(
\left[\begin{array}{c}{\bf Y}_1\\ {\bf Y}_2\end{array}\right]
\right)\le R\}\,.
\ee

\noindent and the problem
\be
\label{prob-W359}
\min_{{\cal Y}\in \wt{W}_R(I,J,2)} \|{\cal Z}-{\cal Y}\|^2_F\,.
\ee

\noindent Since $\wt{W}_R(I,J,2)$ is closed, problem (\ref{prob-W359}) is
guaranteed to have an optimal solution.
We have $S_R(I,J,2)\subset\wt{W}_R(I,J,2)$, see \cite{Ste2}.
Hence, if a best approximation ${\cal X}$ from the set $\wt{W}_R(I,J,2)$
has rank at most $R$, then ${\cal Z}$ has a best rank-$R$ approximation.
Note that problem~(\ref{prob-W359}) is equivalent to finding a best multilinear
rank-$(R,R,2)$ approximation of ${\cal Z}$, with no transformation in the third mode.
Algorithms to solve this problem can be found in \cite{Sav} \cite{Ish}.

Next, we present an algorithm to solve problem (\ref{prob-W359}) by using Givens rotations.
We make use of this algorithm in our proof for cases 3, 5, and 9.
Let ${\cal Y}\in\wt{W}_R(I,J,2)$ have the following SVDs of its matrix unfoldings:
\be
[{\bf Y}_1\,|\,{\bf Y}_2]={\bf U}_1\,{\bf S}_1\,{\bf V}_1^T\,,\quad\quad
\left[\begin{array}{c}{\bf Y}_1\\ {\bf Y}_2\end{array}\right]=
{\bf V}_2\,{\bf S}_2\,{\bf U}_2^T\,,
\ee

\noindent with ${\bf U}_1$ ($I\times I$), ${\bf V}_1$ ($2J\times 2J$),
${\bf V}_2$ ($2I\times 2I$), and ${\bf U}_2$ ($J\times J$) orthonormal.
Since both unfoldings have rank at most $R$, only the first $R$ diagonal
entries of ${\bf S}_1$ and ${\bf S}_2$ are nonzero. It follows that
\be
{\bf U}_1^T\,{\bf Y}_k\,{\bf U}_2=\left[\begin{array}{cc}
{\bf G}_k & {\bf O}_{R,J-R} \\
{\bf O}_{I-R,R} & {\bf O}_{I-R,J-R}\end{array}\right]\,,\quad\quad
k=1,2\,,
\ee

\noindent where ${\bf G}_k$ is $R\times R$, $k=1,2$. Hence, ${\cal
Y}\in\wt{W}_R(I,J,2)$ satisfies ${\bf Y}_k={\bf U}_{1,R}\,{\bf G}_k\,{\bf
U}_{2,R}^T$, $k=1,2$, where ${\bf U}_{1,R}$ ($I\times R$) and ${\bf
U}_{2,R}$ ($J\times R$) consist of the first $R$ columns of ${\bf U}_1$
and ${\bf U}_2$, respectively. Analogous to the GSD algorithm discussed in
section 5.1, problem (\ref{prob-W359}) can be solved by finding orthonormal
${\bf U}_1$ and ${\bf U}_2$ that maximize the Frobenius norm of the first
$R$ rows and columns of ${\bf U}_1^T\,{\bf Z}_k\,{\bf U}_2$, $k=1,2$. The
best approximation ${\cal X}$ from $\wt{W}_R(I,J,2)$ then has slices ${\bf
X}_k={\bf U}_{1,R}\,{\bf G}_k\,{\bf U}_{2,R}^T$, $k=1,2$, where ${\bf
U}_{1,R}$ and ${\bf U}_{2,R}$ consist of the first $R$ columns of ${\bf
U}_1$ and ${\bf U}_2$, respectively, and ${\bf G}_k$ is taken as the first
$R$ rows and columns of ${\bf U}_1^T\,{\bf Z}_k\,{\bf U}_2$, $k=1,2$.

Finding ${\bf U}_1$ and ${\bf U}_2$ can be done by Givens rotations as follows.
We write
\be
\label{eq-alg1opt}
{\bf U}_1^T\,{\bf Z}_k\,{\bf U}_2=\left[\begin{array}{cc}
{\bf G}_k & {\bf L}_k \\
{\bf H}_k & {\bf M}_k\end{array}\right]\,,
\quad\quad k=1,2\,.
\ee

\noindent To maximize $\|{\bf G}_1\|_F^2+\|{\bf G}_2\|_F^2$, iterative Givens rotations
can be used for each pair of rows $(i,j)$, $1\le i\le R$, $R+1\le j\le I$, in
$\left[\begin{array}{c|c}
{\bf G}_1 & {\bf G}_2\\
{\bf H}_1 & {\bf H}_2\end{array}\right]$, and for each pair
of columns $(i,j)$, $1\le i\le R$, $R+1\le j\le J$, in
$\left[\begin{array}{cc}
{\bf G}_1 & {\bf L}_1\\
\hline
{\bf G}_2 & {\bf L}_2\end{array}\right]$. Analogous to the
derivation of first-order optimality conditions for the GSD algorithm in
section 5.1, we obtain the following conditions for stationary points of
problem~(\ref{prob-W359}) in terms of (\ref{eq-alg1opt}):

\begin{itemize}
\item[$\bullet$] All rows of $[{\bf H}_1\,|\,{\bf H}_2]$ are orthogonal to
all rows of $[{\bf G}_1\,|\,{\bf G}_2]$.
\item[$\bullet$] All columns of $\left[\begin{array}{c}
{\bf L}_1\\ {\bf L}_2\end{array}\right]$ are orthogonal to all columns of
$\left[\begin{array}{c} {\bf G}_1\\ {\bf G}_2\end{array}\right]$.
\end{itemize}

\noindent The rank of a best approximation ${\cal X}$ from the set
$\wt{W}_R(I,J,2)$ is equal to the rank of the $R\times R\times 2$ array
${\cal G}$ with slices ${\bf G}_1$ and ${\bf G}_2$. As stated above, rank$({\cal
X})=$ rank$({\cal G})\le R$ implies that ${\cal Z}$ has a best rank-$R$
approximation. The rank of ${\cal G}$ can be checked by making use of
Lemma~\ref{lem-rank}.

We have the following result for the case where rank$({\cal X})=$
rank$({\cal G})>R$.

\begin{thm}
\label{t-359}
Let ${\cal Z}\in\R^{I\times J\times 2}$ with $2\le R<\min(I,J)$ be generic.
Let $({\bf Q}_a,{\bf Q}_b,{\bf R}_1,{\bf R}_2)$ be an optimal solution of the
GSD problem $(\ref{eq-GSD})$ with strictly positive second-order conditions $(\ref{eq-second1})$-$(\ref{eq-second2})$ for $1\le i<j\le R$. If all optimal solutions ${\cal X}$ of
problem $(\ref{prob-W359})$ have {\rm rank}$({\cal X})>R$, then the rank of
the $R\times R\times 2$ array ${\cal R}$ with slices ${\bf R}_1$ and ${\bf
R}_2$ is larger than $R$.\end{thm}

\noindent {\bf Proof.} See Appendix D.\ep

\noindent Analogous to Conjecture~\ref{conjec}, we assume that the second-order conditions (\ref{eq-second1})-(\ref{eq-second2}), $1\le i<j\le R$, are strictly positive for generic ${\cal Z}\in\R^{I\times J\times 2}$ with $2\le R<\min(I,J)$. As in Appendix A, we have found no counterexamples in numerical experiments (results not reported). Analogous to Corollary~\ref{c-8} following from Theorem~\ref{t-8}, we obtain the following.

\begin{cor}
\label{c-359}
Let ${\cal Z}\in\R^{I\times J\times 2}$ with $2\le R<\min(I,J)$ be generic.
If all optimal solutions ${\cal X}$ of problem $(\ref{prob-W359})$ have
{\rm rank}$({\cal X})>R$, then ${\cal Z}$ does not have a best rank-$R$
approximation. \ep \end{cor}

\noindent Corollary~\ref{c-359} implies that we now have an easy-to-check
criterion to determine whether ${\cal Z}$ has a best rank-$R$ approximation
or not. First, compute a best approximation ${\cal X}$ from the set
$\wt{W}_R(I,J,2)$ by using the algorithm with Givens rotations or the algorithms 
in \cite{Sav} \cite{Ish}. A number of runs with random
starting values can be executed to make sure the global maximum is obtained
and the optimal solution ${\cal X}$ is unique. As in case 8, array
${\cal G}$ (corresponding to ${\cal X}$) may be considered a generic
$R\times R\times 2$ array. Hence, rank$({\cal X})=$ rank$({\cal G})$ equals
$R$ or $R+1$, both on sets of positive Lebesgue measure \cite{TBK}. Next,
compute the eigenvalues of ${\bf G}_2{\bf G}_1^{-1}$, which are distinct.
If all eigenvalues are real, then rank$({\cal X})=$ rank$({\cal G})=R$
(Lemma~\ref{lem-rank}) and ${\cal Z}$ has a best rank-$R$ approximation,
which can be taken equal to ${\cal X}$. If some eigenvalues are complex,
then rank$({\cal X})=$ rank$({\cal G})=R+1$ (Lemma~\ref{lem-rank}) and
${\cal Z}$ does not have a best rank-$R$ approximation. Since both
situations occur on sets of positive Lebesgue measure, this completes the
proof of cases 3, 5, and 9 of Table~\ref{tab-1}.

\section{Discussion}
\setcounter{equation}{0}
Using the Generalized Schur Decomposition (GSD) and its relation to the CPD,
we have proven all conjectures of \cite{Ste2} on the (non)existence of
best rank-$R$ approximations for generic $I\times J\times 2$ arrays. Our
main result is that generic $I\times I\times 2$ arrays of rank $I+1$ do not
have a best rank-$I$ approximation. So far, this was only proven for $I=2$
\cite{DSL}, which was the only result on (non)existence of low-rank
approximations for generic three-way arrays in the literature.

In cases 3, 5, 8, and 9 of Table~\ref{tab-1}, existence of a best rank-$R$
approximation holds on a set of positive volume only. For such arrays we
have obtained easy-to-check necessary and sufficient conditions for the
existence of a best rank-$R$ approximation. In case 8, it suffices to solve
problem (\ref{prob-W}) by computing a truncated SVD and computing the
eigenvalues of a corresponding matrix. In cases 3, 5, and 9, problem
(\ref{prob-W359}) needs to be solved, and the
eigenvalues of a matrix corresponding to the optimal solution of
(\ref{prob-W359}) need to be computed. To the author's knowledge, this
is the first time such conditions are formulated.

In our proofs, we have made use of the fact that the GSD solution set
$P_R(I,J,2)$ is the closure of the set $S_R(I,J,2)$ of arrays with rank at
most $R$ \cite{Ste-arxiv}. Unfortunately, this result does not generalize to
$I\times J\times K$ arrays with $K\ge 3$ and the Simultaneous Generalized
Schur Decomposition \cite{LMV}, nor do we know any other closed form
description of the closure of the rank-$R$ set $S_R(I,J,K)$. Hence,
at present the results for $I\times J\times 2$ arrays in this paper do not
seem to be generalizable to $I\times J\times K$ arrays.

\section*{Appendix A: numerical evidence for Conjecture~\ref{conjec}}
\refstepcounter{section}
\setcounter{equation}{0}
\renewcommand{\thesection}{A}
\renewcommand{\theequation}{A.\arabic{equation}}
Evidence for Conjecture~\ref{conjec} is provided by running the GSD algorithm for random 
$I\times I\times 2$ arrays ${\cal Z}$ of rank $I+1$. The GSD algorithm is terminated when the relative decrease in error sum of squares drops below $10^{-9}$. For each array, we record the maximum absolute value (over all pairs $(i,j)$) of the first-order conditions (\ref{eq-stat1}) and (\ref{eq-stat2}), and the minimum value (over all pairs $(i,j)$) of the second-order conditions (\ref{eq-second1}) and (\ref{eq-second2}). According to Conjecture~\ref{conjec}, the latter should be strictly positive for an optimal GSD solution.  To evaluate whether the GSD algorithm has terminated in a local minimum, we also compute the Hessian matrix of second order derivatives. Let the orthonormal ${\bf Q}_a$ and ${\bf Q}_b$ be parameterized by $I(I-1)/2$ Gives rotations each. We consider the Hessian matrix of the problem (\ref{eq-GSD3}), where the variables are the $I(I-1)$ rotation angles. The first-order conditions are then equal to (\ref{eq-stat1}) and (\ref{eq-stat2}), and the second-order conditions (\ref{eq-second1}) and (\ref{eq-second2}) form the diagonal of the Hessian matrix. The Hessian matrix of the GSD problem was not considered in \cite{LMV}. 

Note that the GSD algorithm decreases the objective value (\ref{eq-GSD3}) with every Givens rotation, unless the corresponding second-order condition equals zero (then the objective value remains unchanged; see Lemma~\ref{lem-statorth}). Also, the GSD algorithm is expected to converge to a stationary point, i.e., with zero first-order conditions (\ref{eq-stat1}) and (\ref{eq-stat2}). Indeed, if a first-order condition is nonzero, then a Givens rotation exists that will decrease the objective value. Of course, since the GSD algorithm is terminated after a finite number of iterations, the first-order conditions (\ref{eq-stat1}) and (\ref{eq-stat2}) will not be exactly zero after convergence. 

In \cite{LMV}, the GSD algorithm is initialized by taking ${\bf Q}_a$ and ${\bf Q}_b$ from the generalized real Schur decomposition \cite[theorem 7.7.2]{GVL} computed via QZ iteration, and an extensive simulation study yields no cases of suboptimal GSD solutions. Here, we also consider initialization with random orthonormal ${\bf Q}_a$ and ${\bf Q}_b$ and check whether the GSD algorithm terminates in suboptimal solutions. A GSD solution is suboptimal when a better GSD solution is found for different initial values in the GSD algorithm. For each $I\in\{2,3,4,5,6,7\}$ we generate 10 arrays ${\cal Z}$ and run the GSD algorithm 100 times with random initial values, and 1 time with the QZ initial values. In Table~\ref{tab-sim1} (top half) the results are presented. All runs result in first-order conditions close to zero, with largest absolute value below $10^{-4}$. All runs result in strictly positive second-order conditions, with the smallest value being $0.14$ for $I=3$. It is clear that Conjecture~\ref{conjec} is not violated in these simulations. All runs result in a positive definite Hessian matrix, indicating a local minimum, although for some runs the Hessian has a very small positive eigenvalue. Hence, we do not find any saddle points. \\

\begin{table}[htb]
\begin{center}
\begin{tabular}{|c||c|c|c|c|}
\hline $I$ & suboptimal \% & max(abs(1st)) & min(2nd) & min(eig(Hessian)) \\[2mm]

\hline\hline

$I=2$ & 0 & $7\cdot 10^{-16}$ & 0.15 & 0.0215 \\
$I=3$ & 3.9 & $5\cdot 10^{-5}$ & 0.14 & 0.0082 \\
$I=4$ & 6.4 & $1\cdot 10^{-4}$ & 0.77 & 0.0028 \\
$I=5$ & 15.1 & $8\cdot 10^{-5}$ & 0.59 & 0.0010 \\
$I=6$ & 29.5 & $8\cdot 10^{-5}$ & 0.59 & 0.0017 \\
$I=7$ & 25.2 & $9\cdot 10^{-5}$ & 0.86 & 0.0005 \\

\hline\hline

$I=2$ & - & $9\cdot 10^{-16}$ & 0.01 & 0.0024 \\
$I=3$ & - & $4\cdot 10^{-16}$ & 0.01 & 0.0009 \\
$I=4$ & - & $9\cdot 10^{-5}$ & 0.08 & 0.0002 \\
$I=5$ & - & $1\cdot 10^{-4}$ & 0.29 & $9\cdot 10^{-5}$ \\
$I=6$ & - & $1\cdot 10^{-4}$ & 0.28 & $7\cdot 10^{-6}$ \\
$I=7$ & - & $1\cdot 10^{-4}$ & 0.36 & $1\cdot 10^{-5}$ \\

\hline

\end{tabular}
\end{center}

\caption{Percentage of suboptimal solutions, maximum absolute value of first-order conditions (\ref{eq-stat1})-(\ref{eq-stat2}), minimum value of second-order conditions (\ref{eq-second1})-(\ref{eq-second2}), and smallest eigenvalue of the Hessian for runs of the GSD algorithm for random $I\times I\times 2$ arrays ${\cal Z}$ of rank $I+1$. Top half: 10 arrays per $I$; for each array 101 runs are executed, 100 with random starting values and 1 with QZ starting values; results for optimal runs are reported. Bottom half: 1000 arrays per $I$; for each array 1 run with QZ starting values is executed.}
\label{tab-sim1} \end{table}\vspace{5mm}

\noindent Quite some runs result in a suboptimal local minimum, but the QZ initialized runs result in a suboptimal solution for only 3 out of the 60 arrays. In all following simulations, we use only one QZ initialized run for each array. Next, we generate 1000 arrays ${\cal Z}$ for each $I\in\{2,3,4,5,6,7\}$ and again check the first- and second-order conditions, and the eigenvalues of the Hessian matrix. The results are reported in Table~\ref{tab-sim1} (bottom half). Again the first-order conditions are close to zero, with largest absolute value below $10^{-4}$. The second-order conditions are all positive, but the minimal values for $I=2,3,4$ are rather small. However, they are still $10^{14}$ times larger than the largest first-order condition for $I=2,3$, and  $10^{3}$ times larger for $I=4$. Therefore, we do not consider this a violation of Conjecture~\ref{conjec}. The smallest eigenvalues of the Hessian matrix are positive but small. Again we do not encounter any saddle points. 

For larger values of $I$ the symbolic computation of the Hessian matrix takes a lot of time (for $I=7$ it takes 6 hours on a regular PC). Hence, we omit the computation of the Hessian in the final simulations for $I=2,\ldots,25$. For each value of $I$, we generate 10 arrays ${\cal Z}$ and run the GSD algorithm with QZ initial values. The results are depicted in Figure~\ref{fig-1} below. As can be seen, for all arrays the largest absolute value of the first-order conditions is below $10^{-4}$, and the smallest value of the second-order conditions is strictly positive (with the smallest values being $0.24$ for $I=2$ and $I=3$). Hence, also here Conjecture~\ref{conjec} is not violated.\\

\begin{figure}[h]
\begin{center}
\includegraphics[height=6cm]{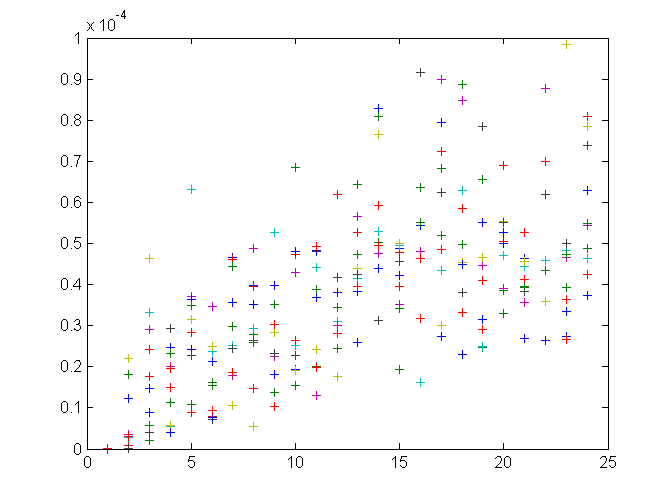}
\includegraphics[height=6cm]{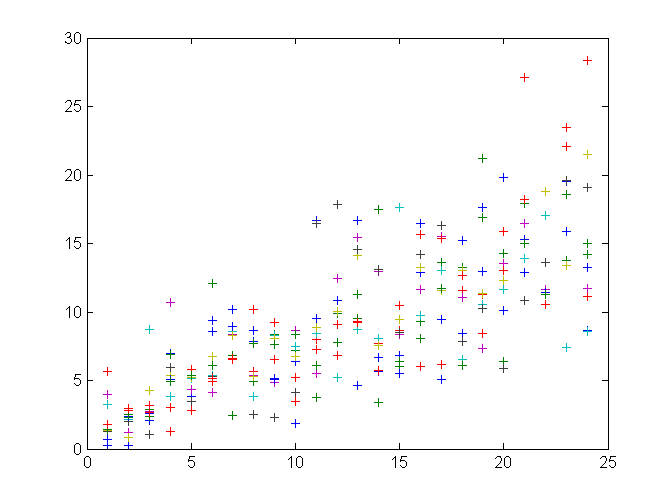}
\caption{Maximum absolute value of first-order conditions (\ref{eq-stat1})-(\ref{eq-stat2}) (left), and minimum value of second-order conditions (\ref{eq-second1})-(\ref{eq-second2}) (right), after running the GSD algorithm for 10 random $I\times I\times 2$ arrays ${\cal Z}$ of rank $I+1$, for $I=2,\ldots,25$.} 
\label{fig-1} \end{center}
\end{figure}

\noindent Note that Conjecture~\ref{conjec} does not imply that zero second-order conditions do not occur for $I\times I\times 2$ arrays ${\cal Z}$ of rank $I+1$. The conjecture only states that this does not occur for {\em generic} $I\times I\times 2$ arrays ${\cal Z}$ of rank $I+1$. Below are two examples with zero second-order conditions for $I=3$. Let
\be
\label{eq-exampleA1}
[\wt{\bf Z}_1\,|\,\wt{\bf Z}_2]=\left[\begin{array}{lcc|lcc}
0 & 1 & 2 & 0 & 2 & -2 \\
0 & 2 & 3 & 0 & 3 & -4 \\
0.1 & 0 & 2 & 0.2 & 0 & -1\end{array}\right]\,.
\ee

\noindent The matrix $\wt{\bf Z}_2\wt{\bf Z}_1^{-1}$ has complex eigenvalues, which implies rank$({\cal Z})=4$ by Lemma~\ref{lem-rank}. It can be verified that (\ref{eq-exampleA1}) satisfies the first-order conditions (\ref{eq-stat1})-(\ref{eq-stat2}) and the second-order conditions (\ref{eq-second1})-(\ref{eq-second2}). However, the second-order condition (\ref{eq-second1}) is zero for $(i,j)=(1,2)$. The Hessian corresponding to (\ref{eq-exampleA1}) has eigenvalues -0.1, 0, 10.0, 22.1, 36.1, and 63.7. Hence, this is an example of a saddle point. 

A second example is obtained by setting the (3,3) entries of $\wt{\bf Z}_1$ and $\wt{\bf Z}_2$ in (\ref{eq-exampleA1}) to zero. Let
\be
\label{eq-exampleA2}
[\wt{\bf Z}_1\,|\,\wt{\bf Z}_2]=\left[\begin{array}{lcc|lcr}
0 & 1 & 2 & 0 & 2 & -2 \\
0 & 2 & 3 & 0 & 3 & -4 \\
0.1 & 0 & 0 & 0.2 & 0 & 0\end{array}\right]\,.
\ee

\noindent Again $\wt{\bf Z}_2\wt{\bf Z}_1^{-1}$ has complex eigenvalues, and the first- and second-order conditions are satisfied. The second-order condition (\ref{eq-second1}) is zero for $(i,j)=(1,2)$, and the second-order condition (\ref{eq-second2}) is zero for $(i,j)=(2,3)$. The corresponding Hessian has eigenvalues 0, 0, 0.01, 20.9, 36.0, and 54.9. Hence, this is not a saddle point. Both examples (\ref{eq-exampleA1}) and (\ref{eq-exampleA2}) are suboptimal GSD solutions and can be improved by using the GSD algorithm with random starting values. 

An alternative way to verify the strict positivity of the second-order conditions is by considering {\em all stationary points} of the GSD problem (\ref{eq-GSD3}), and using a dimensionality argument. Consider the set $Q_I$ of all ${\bf Q}_a$ and ${\bf Q}_b$ satisfying the first-order optimality conditions (\ref{eq-stat1}) and (\ref{eq-stat2}). We rewrite the latter explicitly in terms of  ${\bf Q}_a$ and ${\bf Q}_b$
and obtain
\begin{eqnarray}
\label{eq-setQ}
Q_I=\{{\bf Q}_a,{\bf Q}_b\in\R^{I\times I}&:&{\bf Q}_a^T{\bf Q}_a={\bf Q}_b^T{\bf Q}_b={\bf I}_I\,, \nonumber \\
&& \sum_{r=i}^{j-1} \sum_{k=1}^2 ({\bf q}_{a,i}^T{\bf Z}_k{\bf q}_{b,r})
({\bf q}_{a,j}^T{\bf Z}_k{\bf q}_{b,r})=0\,,\quad 1\le i<j\le I\,, \nonumber \\
&& \sum_{r=i+1}^{j} \sum_{k=1}^2 ({\bf q}_{a,r}^T{\bf Z}_k{\bf q}_{b,i})
({\bf q}_{a,r}^T{\bf Z}_k{\bf q}_{b,j})=0\,,\quad 1\le i<j\le I\,\}\,,
\end{eqnarray}

\noindent where ${\bf q}_{a,r}$ and ${\bf q}_{b,r}$ denote the $r$th columns of ${\bf Q}_a$ and ${\bf Q}_b$, respectively. The number of equations in ${\bf Q}_a^T{\bf Q}_a={\bf Q}_b^T{\bf Q}_b={\bf I}_I$ equals $2\,I(I+1)/2=I(I+1)$. The first-order optimality conditions form $2\,I(I-1)/2=I(I-1)$ equations. Hence, the total number of equations in (\ref{eq-setQ}) equals $I(I+1)+I(I-1)=2\,I^2$, which is equal to the number of unknowns in ${\bf Q}_a$ and ${\bf Q}_b$. If the dimensionality of the set $Q_I$ equals zero, then it contains at most a finite number of ${\bf Q}_a$ and ${\bf Q}_b$. Since the GSD problem (\ref{eq-GSD3}) always has an optimal solution, the set $Q_I$ is not empty. Moreover, due to permutation ambiguities in the GSD, the set $Q_I$ always contains more than one pair of ${\bf Q}_a$ and ${\bf Q}_b$ \cite{LMV} \cite{SDL}. Lemma~\ref{lem-statorth} shows that a second-order condition being zero at a stationary point implies infinitely many optimal ${\bf Q}_a$ or ${\bf Q}_b$ (a one-dimensional set, in fact). Hence, this would contradict the dimensionality of the set $Q_I$ being zero. With the use of a computer algebra package, the number of (real) solutions to the equations defining $Q_I$ could be computed, for some random $I\times I\times 2$ array ${\cal Z}$ of rank $I+1$. When this number is finite, the second-order conditions are strictly positive. However, although this approach was numerically feasible for $I=2$, for $I=3$ no solution to the equations was found after 6 hours on a regular PC using Mathematica 9. Note that for $I=3$ we have 18 4th degree poynomial equations in 18 variables. Another approach to verify that dim$(Q_I)=0$ is by using an algebraic geometry program to compute the dimension of the ideal generated by the polynomial equations in $Q_I$, with the rational numbers as base field for the coefficients. Hence, here ${\bf Z}_1$ and ${\bf Z}_2$ should be random and rational, which approximates the real case. We have used the Macaulay2 package and for $I=2$ the dimensionality of the ideal was indeed zero. However, for $I=3$ computing the dimensionality of the ideal did not produce an answer after several hours on a regular PC. Since the approach via dim$(Q_I)=0$ was not successful, we verified the strict positivity of the second-order conditions using the GSD algorithm instead.

In \cite{LMV} the GSD algorithm is used as an algorithm for simultaneous matrix diagonalization. When the stationary points in $Q_I$ could be computed easily, then a new fast algorithm for simultaneous matrix diagonalization would have been obtained. Moreover, all local and global minima would be known. This would be a huge improvement with respect to all other exisiting algorithms for simultaneous matrix diagonalization. Therefore, we do not expect that computing all solutions to the equations in $Q_I$ is numerically feasible for $I$ not very small.

\section*{Appendix B: proof of Theorem~\ref{t-1}}
\refstepcounter{section}
\setcounter{equation}{0}
\renewcommand{\thesection}{B}
\renewcommand{\theequation}{B.\arabic{equation}}
First, we show that we may assume without loss of generality that
the optimal ${\bf R}_1$ and ${\bf R}_2$ are nonsingular when the second-order conditions (\ref{eq-second1})-(\ref{eq-second2}) are strictly positive. The optimal ${\bf
R}_k$ equals the upper triangular part of $\wt{\bf Z}_k={\bf Q}_a^T\,{\bf
Z}_k\,{\bf Q}_b$, $k=1,2$. A singular ${\bf R}_k$ implies that it has one
or more diagonal entries equal to zero. The second-order optimality
conditions (\ref{eq-second1})-(\ref{eq-second2}) for rotations $(i,i+1)$
are $\tilde{\bf z}_{(i,i)}^T\tilde{\bf z}_{(i,i)}\ge\tilde{\bf z}_{(i+1,i)}^T\tilde{\bf z}_{(i+1,i)}$, $i=1,\ldots,I-1$, and $\tilde{\bf z}_{(i,i)}^T\tilde{\bf z}_{(i,i)}\ge\tilde{\bf z}_{(i,i-1)}^T\tilde{\bf z}_{(i,i-1)}$, $i=2,\ldots,I$.
Hence, positive second-order conditions imply that ${\bf
R}_1$ and ${\bf R}_2$ do not both have a zero on position $(i,i)$,
$i=1,\ldots,I$. In that case, a nonsingular slicemix ${\bf S}$ ($2\times
2$) exists, such that the mixed slices $\wh{\bf R}_k=s_{k1}\,{\bf
R}_1+s_{k2}\,{\bf R}_2$, $k=1,2$, are nonsingular. The matrix ${\bf S}$ may
be taken orthonormal. Then the GSD solution $({\bf Q}_a,{\bf Q}_b,\wh{\bf
R}_1,\wh{\bf R}_2)$ is an optimal solution of the GSD problem for array
$\wh{\cal Z}$ with mixed slices $\wh{\bf Z}_k=s_{k1}\,{\bf
Z}_1+s_{k2}\,{\bf Z}_2$, $k=1,2$. The optimal $\wh{\bf R}_k$, $k=1,2$, are
nonsingular, rank$({\cal Z})=$ rank$(\wh{\cal Z})$, and rank$({\cal R})=$
rank$(\wh{\cal R})$. Hence, a proof of rank$(\wh{\cal R})>I$ implies
rank$({\cal R})>I$. Therefore, positive second-order conditions imply
that we may assume without loss of generality that the optimal ${\bf R}_1$
and ${\bf R}_2$ are nonsingular.

As rank criterion for ${\cal R}$ we use Lemma~\ref{lem-rank} $(ii)$. We
show that
\begin{itemize}
\item[$({\rm B}i)$] ${\bf R}_2{\bf R}_1^{-1}$ has some identical
eigenvalues, \item[$({\rm B}ii)$] ${\bf R}_2{\bf R}_1^{-1}$ does not have
$I$ linearly independent eigenvectors. \end{itemize}

\noindent First, we prove $({\rm B}i)$ by contradiction. Suppose the
eigenvalues of ${\bf R}_2{\bf R}_1^{-1}$ are distinct. Then rank$({\cal
R})=I$ by Lemma~\ref{lem-rank} $(i)$. The eigenvalues of ${\bf R}_2{\bf
R}_1^{-1}$ are equal to the diagonal entries of $\wt{\bf Z}_2$ divided by
those of $\wt{\bf Z}_1$. Statement $({\rm B}i)$ not holding is equivalent
to none of the vectors $\tilde{\bf z}_{(i,i)}$, $i=1,\ldots,I$, being
proportional. The nonsingularity of ${\bf R}_1$ and ${\bf R}_2$ implies
that vectors $\tilde{\bf z}_{(i,i)}$ do not contain zeros, $i=1,\ldots,I$.
Consider the optimality conditions (\ref{eq-stat1}) and (\ref{eq-stat2})
for rotations $(i,i+1)$, which are: $\tilde{\bf z}_{(i,i)}^T\tilde{\bf
z}_{(i+1,i)}=\tilde{\bf z}_{(i+1,i)}^T\tilde{\bf z}_{(i+1,i+1)}=0$. Since
$\tilde{\bf z}_{(i,i)}$ and $\tilde{\bf z}_{(i+1,i+1)}$ are not
proportional and not zero, this implies $\tilde{\bf z}_{(i+1,i)}={\bf
0}$ for $i=1,\ldots,I-1$. Here, ${\bf 0}$ denotes the zero vector in
$\R^2$.

Using this result, we consider (\ref{eq-stat1}) and (\ref{eq-stat2}) for
rotations $(i,i+2)$. These equations now become (the sum terms vanish):
$\tilde{\bf z}_{(i,i)}^T\tilde{\bf z}_{(i+2,i)}=\tilde{\bf
z}_{(i+2,i)}^T{\bf z}_{(i+2,i+2)}=0$. Since $\tilde{\bf z}_{(i,i)}$
and $\tilde{\bf z}_{(i+2,i+2)}$ are not proportional and not zero, this
implies $\tilde{\bf z}_{(i+2,i)}={\bf 0}$ for $i=1,\ldots,I-2$. By
consecutively considering rotations $(i,i+q)$ in this way, it is clear that
we obtain $\tilde{\bf z}_{(j,i)}={\bf 0}$ for $1\le i<j\le I$. 
This implies that $\wt{\bf Z}_k={\bf
Q}_a^T\,{\bf Z}_k\,{\bf Q}_b$, $k=1,2$, are upper triangular. Moreover,
$\wt{\bf Z}_k={\bf R}_k$, $k=1,2$, and rank$({\cal Z})=$ rank$({\cal
R})=I$, which contradicts the assumption of rank$({\cal Z})=I+1$ in
Theorem~\ref{t-1}. Hence, we have proven that $({\rm B}i)$ holds.

For the proof of $({\rm B}ii)$ we first reorder the eigenvalues of ${\bf
R}_2{\bf R}_1^{-1}$ such that identical eigenvalues appear in contiguous
groups. The reordering can be done within the GSD; see Kressner \cite{Kre}.
The proof of $({\rm B}ii)$ is by contradiction. We suppose that ${\bf
R}_2{\bf R}_1^{-1}$ has $I$ linearly independent eigenvectors, which
implies rank$({\cal R})=I$ by Lemma~\ref{lem-rank} $(i)$. Hence, for
eigenvalue $\lambda$ of ${\bf R}_2{\bf R}_1^{-1}$ with multiplicity $d$
there are $d$ linearly independent eigenvectors. In other words, rank$({\bf
R}_2{\bf R}_1^{-1}-\lambda\,{\bf I}_I)=I-d$, which is equivalent to
rank$({\bf R}_2-\lambda\,{\bf R}_1)=I-d$. The diagonal of ${\bf
R}_2-\lambda\,{\bf R}_1$ consists of $I-d$  nonzero entries and $d$ zeros
which form the diagonal of a $d\times d$ strictly upper triangular block ${\bf V}$.
We consider ${\bf R}_2-\lambda\,{\bf R}_1$ as a block upper triangular
matrix. Its rank is at least equal to the sum of the ranks of its diagonal
blocks. These blocks are $I-d$ nonzero scalars and the block ${\bf V}$.
Hence, the lower bound for the rank of ${\bf R}_2-\lambda\,{\bf R}_1$
includes the value $I-d$ only if rank$({\bf V})=0$. Conversely, row and
column operations can be used to show that rank$({\bf V})=0$ implies
rank$({\bf R}_2-\lambda\,{\bf R}_1)=I-d$. Rank$({\bf V})=0$ is equivalent
to all vectors $\tilde{\bf z}_{(m,n)}$ contained in the upper triangular
part of the $d\times d$ block (with $m\le n$) being either proportional or
zero. Proposition~\ref{p-2} below shows that the vectors $\tilde{\bf
z}_{(m,n)}$ with $m>n$ in the block are zero. To summarize, let the
corresponding $d\times d$ diagonal block of $\wt{\bf Z}_k={\bf Q}_a^T\,{\bf
Z}_k\,{\bf Q}_b$ be denoted by ${\bf W}_k$, $k=1,2$. Then ${\bf W}_1$ and
${\bf W}_2$ are upper triangular and ${\bf W}_2-\lambda\,{\bf W}_1={\bf
O}_{d,d}$, where ${\bf O}_{d,d}$ denotes the $d\times d$ zero matrix.

The last part of the proof of $({\rm B}ii)$ is similar to the proof of
$({\rm B}i)$: we show that $\wt{\bf Z}_k$, $k=1,2$, are upper triangular,
which implies rank$({\cal Z})=$ rank$({\cal R})=I$. The contradiction with
rank$({\cal Z})=I+1$ then implies that
${\bf R}_2{\bf R}_1^{-1}$ does not have $I$ linearly independent
eigenvectors. The diagonal of $\wt{\bf Z}_k$ consists of blocks ${\bf
W}_k^{(1)},\ldots,{\bf W}_k^{(L)}$ corresponding to $L$ distinct
eigenvalues of ${\bf R}_2{\bf R}_1^{-1}$. Let blocks ${\bf W}_k^{(l)}$,
$k=1,2$, have size $d_l\times d_l$, $l=1,\ldots,L$. A $1\times 1$ block
corresponds to a unique eigenvalue, and a $d_l\times d_l$ block with
$d_l\ge 2$ corresponds to an eigenvalue with multiplicity $d_l$. From
Proposition~\ref{p-2} we know that the $d_l\times d_l$ blocks ${\bf
W}_k^{(l)}$, $k=1,2$, are upper triangular. Consider the optimality
conditions (\ref{eq-stat1}) and (\ref{eq-stat2}) for rotations $(i,i+1)$
such that entry $(i,i+1)$ is not part of any block ${\bf W}_k^{(l)}$. These
are: $\tilde{\bf z}_{(i,i)}^T\tilde{\bf z}_{(i+1,i)}=\tilde{\bf
z}_{(i+1,i)}^T\tilde{\bf z}_{(i+1,i+1)}=0$. Since $\tilde{\bf z}_{(i,i)}$
and $\tilde{\bf z}_{(i+1,i+1)}$ are not proportional (they are not in a
block ${\bf W}_k^{(l)}$) and not zero, this implies $\tilde{\bf
z}_{(i+1,i)}={\bf 0}$. Note that $\tilde{\bf z}_{(i+1,i)}={\bf 0}$ for
all $(i+1,i)$ in block ${\bf W}_k^{(l)}$ by Proposition~\ref{p-2}.

As in the proof of $({\rm B}i)$, next we consider (\ref{eq-stat1}) and
(\ref{eq-stat2}) for rotations $(i,i+2)$ such that entry $(i,i+2)$ is not
in a block ${\bf W}_k^{(l)}$. These equations now become (the sum terms
vanish): $\tilde{\bf z}_{(i,i)}^T\tilde{\bf z}_{(i+2,i)}=\tilde{\bf
z}_{(i+2,i)}^T{\bf z}_{(i+2,i+2)}=0$. Since $\tilde{\bf z}_{(i,i)}$ and
$\tilde{\bf z}_{(i+2,i+2)}$ are not proportional and not zero, this
implies $\tilde{\bf z}_{(i+2,i)}={\bf 0}$. Proceeding in the same way, we
obtain $\tilde{\bf z}_{(j,i)}={\bf 0}$ for $1\le i<j\le I$, 
which implies that $\wt{\bf Z}_k$, $k=1,2$, are upper
triangular. This completes the proof of $({\rm B}ii)$.

It remains to state and prove Proposition~\ref{p-2}.

\newpage
\begin{propos}
\label{p-2}
Let ${\cal Z}\in\R^{I\times I\times 2}$ be generic with {\rm rank}$({\cal
Z})=I+1$. Let $({\bf Q}_a,{\bf Q}_b,{\bf R}_1,{\bf R}_2)$ be an optimal
solution of the GSD problem $(\ref{eq-GSD})$, with nonsingular ${\bf R}_1$
and ${\bf R}_2$. For $d\ge 2$, let ${\bf W}_k$ be a $d\times d$ diagonal
block of ${\bf Q}_a^T\,{\bf Z}_k\,{\bf Q}_b$, $k=1,2$, such that the upper
triangular part of ${\bf W}_2-\lambda\,{\bf W}_1$ is zero for some
$\lambda\neq 0$. Then ${\bf W}_1$ and ${\bf W}_2$ are upper triangular.
\end{propos}

\noindent {\bf Proof.} The proof is by induction on $d$. Let ${\bf W}_k^{(h)}$ 
be the $d\times d$ submatrix of $\wt{\bf Z}_k={\bf Q}_a^T\,{\bf Z}_k\,{\bf Q}_b$ consisting of rows
$h,\ldots,h+d-1$ and columns $h,\ldots,h+d-1$. Recall the definition of the vectors $\tilde{\bf
z}_{(m,n)}$ in (\ref{eq-tildez}). Note that since ${\bf R}_k$ are
nonsingular, $k=1,2$, the vectors $\tilde{\bf z}_{(i,i)}$ do not contain
zeros, $i=1,\ldots,I$.

First, we consider $d=2$. We write the $2\times 2\times 2$ array ${\cal W}^{(h)}$ with
unfolding $[{\bf W}_1^{(h)}\;{\bf W}_2^{(h)}]$ in terms of its 3-mode vectors 
$\tilde{\bf z}_{(m,n)}$ as
\be
\label{eq-22}
\left[\begin{array}{cc}
\tilde{\bf z}_{(h,h)} & \tilde{\bf z}_{(h,h+1)}\\
\tilde{\bf z}_{(h+1,h)} & \tilde{\bf z}_{(h+1,h+1)}\end{array}\right]\,,
\ee

\noindent where $\tilde{\bf z}_{(h,h)}$ and $\tilde{\bf z}_{(h+1,h+1)}$ are
proportional, and $\tilde{\bf z}_{(h,h+1)}$ is either zero or proportional
to $\tilde{\bf z}_{(h,h)}$. The proof is by contradiction. 
Suppose $\tilde{\bf z}_{(h+1,h)}\neq {\bf 0}$. 
By the proportionality of $\tilde{\bf z}_{(h,h+1)}$ and $\tilde{\bf z}_{(h,h)}$, 
an orthonormal rotation of the rows of ${\bf W}_k^{(h)}$, $k=1,2$, 
exists that makes $\tilde{\bf z}_{(h,h+1)}$ zero. 
This rotation is not needed when $\tilde{\bf z}_{(h,h+1)}={\bf 0}$. Next, swapping
rows and columns yields upper triangular blocks ${\bf W}_k^{(h)}$, $k=1,2$. 
This implies a better GSD solution has been found, which is a contradiction. 
Hence, it follows that $\tilde{\bf z}_{(h+1,h)}={\bf 0}$. Note that the transformations
used do not affect the GSD objective function outside of the blocks 
${\bf W}_k^{(h)}$, $k=1,2$. This completes the proof for $d=2$.

Next, we consider $d=3$. We write the $3\times 3\times 2$ array ${\cal W}^{(h)}$ with
unfolding $[{\bf W}_1^{(h)}\;{\bf W}_2^{(h)}]$ as 
\be
\label{eq-33}
\left[\begin{array}{ccc}
\tilde{\bf z}_{(h,h)} & \tilde{\bf z}_{(h,h+1)} & \tilde{\bf z}_{(h,h+2)} \\
\tilde{\bf z}_{(h+1,h)} & \tilde{\bf z}_{(h+1,h+1)} & \tilde{\bf z}_{(h+1,h+2)} \\
\tilde{\bf z}_{(h+2,h)} & \tilde{\bf z}_{(h+2,h+1)} & \tilde{\bf z}_{(h+2,h+2)}
\end{array}\right]\,,
\ee

\noindent where $\tilde{\bf z}_{(h+i,h+i)}$, $i=0,1,2$, are proportional, and
$\tilde{\bf z}_{(h+i,h+j)}$ with $0\le i<j\le 2$ are either zero or proportional to
$\tilde{\bf z}_{(h,h)}$. The proof for $d=2$ applies to the subblock
consisting of the first two rows and columns, and to the subblock consisting
of the last two rows and columns. This implies $\tilde{\bf
z}_{(h+1,h)}=\tilde{\bf z}_{(h+2,h+1)}={\bf 0}$. The proof is by contradiction.
Suppose $\tilde{\bf z}_{(h+2,h)}\neq{\bf 0}$.
Let `Row$(i,j)$' to denote an orthonormal rotation of rows $i$ and $j$ of ${\bf W}_k^{(h)}$, $k=1,2$.
Next, we apply the following sequence of orthonormal row rotations:
\bdm
\left[\begin{array}{ccc}
\tilde{\bf z}_{(h,h)} & \tilde{\bf z}_{(h,h+1)} & \tilde{\bf z}_{(h,h+2)} \\
{\bf 0} & \tilde{\bf z}_{(h+1,h+1)} & \tilde{\bf z}_{(h+1,h+2)} \\
\tilde{\bf z}_{(h+2,h)} & {\bf 0} & \tilde{\bf z}_{(h+2,h+2)}
\end{array}\right]\;\stackrel{{\rm Row}(h,h+1)}{\longrightarrow}\;
\left[\begin{array}{ccc}
\bar{\bf z}_{(h,h)} & {\bf 0} & \bar{\bf z}_{(h,h+2)} \\
\bar{\bf z}_{(h+1,h)} & \bar{\bf z}_{(h+1,h+1)} & \bar{\bf z}_{(h+1,h+2)} \\
\tilde{\bf z}_{(h+2,h)} & {\bf 0} & \tilde{\bf z}_{(h+2,h+2)}
\end{array}\right]
\edm
\bdm
\stackrel{{\rm Row}(h,h+2)}{\longrightarrow}\;
\left[\begin{array}{ccc}
\bar{\bar{{\bf z}}}_{(h,h)} & {\bf 0} & {\bf 0} \\
\bar{\bf z}_{(h+1,h)} & \bar{\bf z}_{(h+1,h+1)} & \bar{\bf z}_{(h+1,h+2)} \\
\bar{\bf z}_{(h+2,h)} & {\bf 0} & \bar{\bf z}_{(h+2,h+2)}
\end{array}\right]\,.
\edm

\noindent Hence, first we rotate rows $h$ and $h+1$ such that $\tilde{\bf
z}_{(h,h+1)}$ becomes zero (when $\tilde{\bf z}_{(h,h+1)}$ is not zero
already). Note that $\bar{\bf z}_{(h,h+2)}$ and $\tilde{\bf z}_{(h+2,h+2)}$ are
proportional or $\bar{\bf z}_{(h,h+2)}={\bf 0}$. Then we rotate rows $h$ and $h+2$
to make $\bar{\bf z}_{(h,h+2)}={\bf 0}$. After this, swapping columns $h+1$ and $h+2$ and
swapping rows $h+1$ and $h+2$ makes the blocks lower triangular. Then reversing the
order of the rows and reversing the order of the columns makes the blocks
upper triangular. This yields a better GSD solution, which is a contradiction. 
Hence, we obtain that $\tilde{\bf z}_{(h+2,h)}={\bf 0}$. Note that none of the transformations
used affects the GSD objective function outside of the blocks ${\bf W}_k^{(h)}$, $k=1,2$.
This completes the proof for $d=3$.

Next, we assume the result holds for $d$ and prove it for $d+1$. By the induction
hypothesis, the $(d+1)\times (d+1)\times 2$ array ${\cal W}^{(h)}$ is of the form
\be
\left[\begin{array}{ccccc}
\tilde{\bf z}_{(h,h)} & \tilde{\bf z}_{(h,h+1)} & \ldots & \ldots &
\tilde{\bf z}_{(h,h+d)} \\
{\bf 0} & \tilde{\bf z}_{(h+1,h+1)} & & & \vdots \\
\vdots & \ddots & \ddots & & \vdots \\
{\bf 0} & {\bf 0} & \ddots & \ddots & \vdots \\
\tilde{\bf z}_{(h+d,h)} & {\bf 0} & \ldots & {\bf 0} &
\tilde{\bf z}_{(h+d,h+d)}
\end{array}\right]\,.
\ee

\noindent The proof is by contradiction. Suppose $\tilde{\bf z}_{(h+d,h)}\neq{\bf 0}$.
Next, we apply consecutive rotations of rows $h$ and $h+i$ of ${\bf W}_k^{(h)}$, $k=1,2$, 
to make the current $\tilde{\bf z}_{(h,h+i)}$ zero (if it is not zero already),
for $i=1,\ldots,d$. This yields the form
\be
\left[\begin{array}{ccccc}
\bar{\bf z}_{(h,h)} & {\bf 0} & \ldots & \ldots & {\bf 0} \\
\bar{\bf z}_{(h+1,h)} & \bar{\bf z}_{(h+1,h+1)} & \ldots & \ldots
& \bar{\bf z}_{(h+1,h+d)}\\
\vdots & {\bf 0} & \ddots & & \vdots \\
\vdots & \vdots & \ddots & \ddots & \vdots \\
\bar{\bf z}_{(h+d,h)} & {\bf 0} & \ldots & {\bf 0} &
\bar{\bf z}_{(h+d,h+d)}
\end{array}\right]\,.
\ee

\noindent Reversing the order of columns $h+1,\ldots,h+d$ and reversing the
order of rows $h+1,\ldots,h+d$ then makes the blocks lower triangular.
Next, reversing the order of the rows and reversing the order of the
columns makes the blocks upper triangular. This yields a better GSD
solution, which is a contradiction. Hence, we obtain 
$\tilde{\bf z}_{(h+d,h)}={\bf 0}$. Again, note that none of the transformations
used affects the GSD objective function outside of the blocks ${\bf W}_k^{(h)}$, $k=1,2$.
This completes the proof of Proposition~\ref{p-2}.
\ep

\section*{Appendix C: proof of Theorem~\ref{t-8}}
\refstepcounter{section}
\setcounter{equation}{0}
\renewcommand{\thesection}{C}
\renewcommand{\theequation}{C.\arabic{equation}}
The structure of the proof is analogous to the proof of Theorem~\ref{t-1}.
As explained in Appendix B, the strictly positive second-order optimality conditions
(\ref{eq-second1})-(\ref{eq-second2}) for $1\le i<j\le R$ imply that we may assume
without loss of generality that ${\bf R}_1$ and ${\bf R}_2$ of an optimal
GSD solution are nonsingular. As rank criterion for ${\cal R}$ we use
again Lemma~\ref{lem-rank} $(ii)$. We show that
\begin{itemize}
\item[$({\rm C}i)$] ${\bf R}_2{\bf R}_1^{-1}$ has some identical
eigenvalues, \item[$({\rm C}ii)$] ${\bf R}_2{\bf R}_1^{-1}$ does not have
$R$ linearly independent eigenvectors. \end{itemize}

\noindent We write
\be
\label{eq-Zt8}
\wt{\bf Z}_k=\wt{\bf Q}_a^T\,{\bf Z}_k\,{\bf Q}_b=\left[
\begin{array}{c}
\wt{\bf G}_k\\
\wt{\bf H}_k\end{array}\right]\,,\quad\quad k=1,2\,,
\ee

\noindent where $\wt{\bf G}_k$ is $R\times R$, and $\wt{\bf H}_k$ is
$(I-R)\times R$, $k=1,2$. The GSD algorithm finds $\wt{\bf Q}_a$ and
${\bf Q}_b$ such that the Frobenius norm of the upper triangular parts
of $\wt{\bf G}_k$, $k=1,2$, is maximized. The optimal ${\bf R}_k$ is taken
as the upper triangular part of $\wt{\bf G}_k$, $k=1,2$. Note that $\wt{\bf
Q}_b={\bf Q}_b$ since $R=J$.

First, we prove $({\rm C}i)$ by contradiction. Suppose $({\rm C}i)$ does
not hold, i.e., all eigenvalues of ${\bf R}_2{\bf R}_1^{-1}$ are distinct,
which implies rank$({\cal R})=R$ by Lemma~\ref{lem-rank} $(i)$.
As in the proof of Theorem~\ref{t-1}, optimality conditions
(\ref{eq-stat12})--(\ref{eq-stat22}) then imply that $\wt{\bf G}_k$ are
upper triangular, $k=1,2$. We write $\wt{\bf G}_k={\bf R}_k$, $k=1,2$. The
SVD of $[{\bf Z}_1\,|\,{\bf Z}_2]$ is given as ${\bf U}\,{\bf S}\,{\bf
V}^T$, with ${\bf U}^T{\bf U}={\bf I}_I$ and ${\bf V}^T{\bf V}={\bf
I}_{2J}$. The best approximation ${\cal X}$ from the set
$W_R(I,J,2)$ is given by the truncated SVD and has slices $[{\bf
X}_1\,|\,{\bf X}_2]={\bf U}_R\,{\bf S}_R\,{\bf V}_R^T={\bf U}_R\,[{\bf
S}_R{\bf V}_{R,1}^T\,|\,{\bf S}_R{\bf V}_{R,2}^T]$. The rank of ${\cal X}$
is equal to the rank of the $R\times R\times 2$ array ${\cal V}_R$ with
$R\times R$ slices ${\bf S}_R{\bf V}_{R,1}^T$ and ${\bf S}_R{\bf
V}_{R,2}^T$. In Theorem~\ref{t-8} it is assumed that rank$({\cal X})=$
rank$({\cal V}_R)>R$.

Note that $\wt{\bf Z}_k= \wt{\bf Q}_a^T\,{\bf Z}_k\,{\bf Q}_b=({\bf
U}^T\wt{\bf Q}_a)^T\,({\bf U}^T{\bf Z}_k)\,{\bf Q}_b$, $k=1,2$. We write
\be
\label{eq-Zt8mat}
[\wt{\bf Z}_1\,|\,\wt{\bf Z}_2]=({\bf U}^T\wt{\bf Q}_a)^T\,({\bf SV}^T)\,
\left[\begin{array}{cc}
{\bf Q}_b & {\bf O}_{J,J}\\
{\bf O}_{J,J} & {\bf Q}_b\end{array}\right]=
\left[\begin{array}{cc}
{\bf R}_1 & {\bf R}_2\\
\wt{\bf H}_1 & \wt{\bf H}_2\end{array}\right]\,.
\ee

\noindent The optimality conditions (\ref{eq-stat3}), together with ${\bf
R}_1$ and ${\bf R}_2$ being upper triangular, yield that all rows of
$[\wt{\bf H}_1\;\wt{\bf H}_2]$ are orthognal to all rows of $[{\bf
R}_1\;{\bf R}_2]$. Postmultiplying (\ref{eq-Zt8mat}) by its transpose yields
\be
\label{eq-eigen8}
({\bf U}^T\wt{\bf Q}_a)^T\,{\bf SS}^T\,({\bf U}^T\wt{\bf Q}_a)=
\left[\begin{array}{cc}
{\bf P}_1 & {\bf O}_{R,I-R}\\
{\bf O}_{I-R,R} & {\bf P}_2\end{array}\right]\,,
\ee

\noindent where ${\bf P}_1$ is $R\times R$ and symmetric and nonsingular
(since $[{\bf R}_1\;{\bf R}_2]$ has full row rank due to nonsingularity of
${\bf R}_k$), and ${\bf P}_2$ is $(I-R)\times (I-R)$ and symmetric. Matrix
${\bf S}$ is $I\times 2J$ and contains the $\min(I,2J)$ nonzero singular
values on its diagonal. Note that all singular values are nonzero and distinct since
$[{\bf Z}_1\,|\,{\bf Z}_2]$ is generic.

First, suppose $I\le 2J$. Then ${\bf SS}^T$ is $I\times I$ diagonal and
nonsingular. Since the diagonal entries of ${\bf SS}^T$ are positive and distinct, and
${\bf U}^T\wt{\bf Q}_a$ is orthonormal, equation (\ref{eq-eigen8}) is the
eigendecomposition of an $I\times I$ symmetric matrix that is nonsingular.
In fact, (\ref{eq-eigen8}) can be obtained as the superposition of the
eigendecompositions of ${\bf P}_1$ and ${\bf P}_2$ (both nonsingular).
This implies that
\be
\label{eq-eigenvecsuper8}
{\bf U}^T\wt{\bf Q}_a=\left[\begin{array}{cc}
\wt{\bf Q}_a^{(1)} & {\bf O}_{R,I-R}\\
{\bf O}_{I-R,R} & \wt{\bf Q}_a^{(2)}\end{array}\right]\,,
\ee

\noindent where $\wt{\bf Q}_a^{(1)}$ ($R\times R$) and $\wt{\bf Q}_a^{(2)}$
($(I-R)\times (I-R)$) are orthonormal. From (\ref{eq-Zt8mat}) we then obtain
\begin{eqnarray}
{\bf S}\,{\bf V}^T &=& \left[\begin{array}{cc}
\wt{\bf Q}_a^{(1)} & {\bf O}_{R,I-R}\\
{\bf O}_{I-R,R} & \wt{\bf Q}_a^{(2)}\end{array}\right]\;
\left[\begin{array}{cc}
{\bf R}_1 & {\bf R}_2\\
\wt{\bf H}_1 & \wt{\bf H}_2\end{array}\right]\;
\left[\begin{array}{cc}
{\bf Q}_b^T & {\bf O}_{J,J}\\
{\bf O}_{J,J} & {\bf Q}_b^T\end{array}\right] \nonumber \\[2mm]
&=& \left[\begin{array}{cc}
\wt{\bf Q}_a^{(1)}{\bf R}_1{\bf Q}_b^T & \wt{\bf Q}_a^{(1)}{\bf R}_2{\bf
Q}_b^T \\
\wt{\bf Q}_a^{(2)}\wt{\bf H}_1{\bf Q}_b^T & \wt{\bf Q}_a^{(2)}\wt{\bf
H}_2{\bf Q}_b^T\end{array}\right]\,.
\end{eqnarray}

\noindent This implies that array ${\cal V}_R$ has slices ${\bf S}_R{\bf
V}_{R,k}^T=\wt{\bf Q}_a^{(1)}\,{\bf R}_k\,{\bf Q}_b^T$, $k=1,2$. Hence,
rank$({\cal V}_R)=$ rank$({\cal R})=R$, which contradicts
rank$({\cal X})=$ rank$({\cal V}_R)>R$. Hence, $\wt{\bf G}_k$, $k=1,2$,
cannot be upper triangular and $({\rm C}i)$ must hold.

Next, suppose $I>2J$. Then ${\bf SS}^T$ is $I\times I$ diagonal with the
first $2J$ diagonal entries positive and distinct and the last $I-2J$ diagonal entries
zero. In (\ref{eq-eigen8}), matrix ${\bf P}_1$ is nonsingular and ${\bf
P}_2$ has rank $J=R$. As above, it follows that
\be
\label{eq-eigenvec8}
({\bf U}^T\wt{\bf Q}_a)^T=\left[\begin{array}{ccc}
(\wt{\bf Q}_a^{(1)})^T & {\bf O}_{R,J} & {\bf N}_1 \\
{\bf O}_{I-R,R} & (\wh{\bf Q}_a^{(2)})^T & {\bf N}_2 \end{array}\right]\,,
\ee

\noindent where $(\wh{\bf Q}_a^{(2)})^T$ is $(I-R)\times J$ and contains
the eigenvectors of ${\bf P}_2$, and ${\bf N}_1$ and ${\bf N}_2$ have
$I-2J$ columns. Since the matrix ${\bf U}^T\wt{\bf Q}_a$ is orthonormal,
and $\wt{\bf Q}_a^{(1)}$ is orthonormal, it follows that ${\bf N}_1$ is
zero. Hence, ${\bf U}^T\wt{\bf Q}_a$ is of the same form as in
(\ref{eq-eigenvecsuper8}), and the remaining part of the proof is as above.
This completes the proof of $({\rm C}i)$.

Finally, we prove $({\rm C}ii)$ by contradiction. Suppose ${\bf R}_2{\bf
R}_1^{-1}$ has $R$ linearly independent eigenvectors. Hence, rank$({\cal
R})=R$ by Lemma~\ref{lem-rank} $(i)$. As in the proof of Theorem~\ref{t-1},
the optimality conditions (\ref{eq-stat12})--(\ref{eq-stat22}) then imply
that $\wt{\bf G}_k$, $k=1,2$, are upper triangular. This yields the
contradiction rank$({\cal X})=$ rank$({\cal V}_R)=$ rank$({\cal R})=R$ as
shown in the proof of $({\rm C}i)$ above. Hence, $({\rm C}ii)$ must hold.
This completes the proof of Theorem~\ref{t-8}.

\section*{Appendix D: proof of Theorem~\ref{t-359}}
\refstepcounter{section}
\setcounter{equation}{0}
\renewcommand{\thesection}{D}
\renewcommand{\theequation}{D.\arabic{equation}}
The structure of the proof is analogous to the proof of Theorem~\ref{t-8}.
As before, the strictly positive second-order conditions imply that 
we may assume without loss of generality that ${\bf R}_1$ and
${\bf R}_2$ of an optimal GSD solution are nonsingular. As rank criterion
for ${\cal R}$ we use again Lemma~\ref{lem-rank} $(ii)$. We show that
\begin{itemize}
\item[$({\rm D}i)$] ${\bf R}_2{\bf R}_1^{-1}$ has some identical
eigenvalues, \item[$({\rm D}ii)$] ${\bf R}_2{\bf R}_1^{-1}$ does not have
$R$ linearly independent eigenvectors. \end{itemize}

\noindent We write
\be
\label{eq-Zt359}
\wt{\bf Z}_k=\wt{\bf Q}_a^T\,{\bf Z}_k\,\wt{\bf Q}_b=
\left[\begin{array}{cc}
\wt{\bf G}_k & \wt{\bf L}_k \\
\wt{\bf H}_k & \wt{\bf M}_k\end{array}\right]\,,
\quad\quad k=1,2\,.
\ee

\noindent where $\wt{\bf G}_k$ is $R\times R$, $\wt{\bf H}_k$ is
$(I-R)\times R$, $\wt{\bf L}_k$ is $R\times (J-R)$, and $\wt{\bf M}_k$ is
$(I-R)\times (J-R)$, $k=1,2$. The GSD algorithm finds $\wt{\bf Q}_a$ and
$\wt{\bf Q}_b$ such that the Frobenius norm of the upper triangular parts
of $\wt{\bf G}_k$, $k=1,2$, is maximized. The optimal ${\bf R}_k$ is taken
as the upper triangular part of $\wt{\bf G}_k$, $k=1,2$. The optimal ${\bf
Q}_a$ and ${\bf Q}_b$ of the GSD solution consist of the first $R$ columns
of $\wt{\bf Q}_a$ and $\wt{\bf Q}_b$, respectively.

First, we prove $({\rm D}i)$ by contradiction. Suppose all eigenvalues of
${\bf R}_2{\bf R}_1^{-1}$ are distinct, which implies rank$({\cal R})=R$ by
Lemma~\ref{lem-rank} $(i)$. As in the proof of Theorem~\ref{t-1}, optimality
conditions (\ref{eq-stat12})--(\ref{eq-stat22}) then imply that $\wt{\bf
G}_k$ are upper triangular, $k=1,2$. We write $\wt{\bf G}_k={\bf R}_k$,
$k=1,2$. The optimality conditions (\ref{eq-stat3})--(\ref{eq-stat4}),
together with the upper triangularity of ${\bf R}_1$ and ${\bf R}_2$, yield
that:
\begin{itemize}
\item[$\bullet$] All rows of $[\wt{\bf H}_1\,|\,\wt{\bf H}_2]$ are
orthogonal to all rows of $[{\bf R}_1\,|\,{\bf R}_2]$. \item[$\bullet$] All
columns of $\left[\begin{array}{c} \wt{\bf L}_1\\ \wt{\bf
L}_2\end{array}\right]$ are orthogonal to all columns of
$\left[\begin{array}{c} {\bf R}_1\\ {\bf R}_2\end{array}\right]$.
\end{itemize}

\noindent Note that these are also the conditions for a stationary point of
problem (\ref{prob-W359}). We have
\be
\label{eq-unfold359-1}
\wt{\bf Q}_a^T\,[{\bf Z}_1{\bf Q}_b|\,{\bf Z}_2{\bf Q}_b]=
\left[\begin{array}{cc}
{\bf R}_1 & {\bf R}_2\\
\wt{\bf H}_1 & \wt{\bf H}_2\end{array}\right]\,.
\ee

\noindent Let the SVD of $[{\bf Z}_1{\bf Q}_b|\,{\bf Z}_2{\bf Q}_b]$
($I\times 2R$) be given by ${\bf V}_1\,{\bf S}_1\,{\bf W}_1^T$, with ${\bf
V}_1^T{\bf V}_1={\bf I}_I$ and ${\bf W}_1^T{\bf W}_1={\bf I}_{2R}$.
Postmultiplying (\ref{eq-unfold359-1}) by its transpose yields
\be
\label{eq-eigen359-1}
({\bf V}_1^T\wt{\bf Q}_a)^T\,{\bf S}_1{\bf S}_1^T\,({\bf V}_1^T\wt{\bf
Q}_a)= \left[\begin{array}{cc} {\bf P}_1 & {\bf O}_{R,I-R}\\ {\bf
O}_{I-R,R} & {\bf P}_2\end{array}\right]\,, \ee

\noindent where ${\bf P}_1$ is $R\times R$ and symmetric and nonsingular
(since $[{\bf R}_1\;{\bf R}_2]$ has full row rank due to nonsingularity
of ${\bf R}_k$), and ${\bf P}_2$ is $(I-R)\times (I-R)$ and symmetric.
Matrix ${\bf S}_1$ is $I\times 2R$ and contains the $\min(I,2R)$ nonzero
singular values on its diagonal. Note that all singular values
are nonzero and distinct since $[{\bf Z}_1\,|\,{\bf Z}_2]$ is generic and ${\bf Q}_b$
has full column rank $R$.

Suppose $I\le 2R$. Then ${\bf S}_1{\bf S}_1^T$ is $I\times I$
diagonal and nonsingular. Since the diagonal entries of ${\bf S}_1{\bf
S}_1^T$ are positive and distinct, and ${\bf V}_1^T\wt{\bf Q}_a$ is orthonormal,
equation (\ref{eq-eigen359-1}) is the eigendecomposition of an $I\times I$
symmetric matrix that is nonsingular. In fact, (\ref{eq-eigen359-1}) can be
obtained as the superposition of the eigendecompositions of ${\bf P}_1$ and
${\bf P}_2$ (both nonsingular). This implies that
\be
\label{eq-eigenvecsuper359-1}
{\bf V}_1^T\wt{\bf Q}_a=\left[\begin{array}{cc}
\wt{\bf Q}_a^{(1)} & {\bf O}_{R,I-R}\\
{\bf O}_{I-R,R} & \wt{\bf Q}_a^{(2)}\end{array}\right]\,,
\ee

\noindent where $\wt{\bf Q}_a^{(1)}$ ($R\times R$) and $\wt{\bf Q}_a^{(2)}$
($(I-R)\times (I-R)$) are orthonormal. Hence, $\wt{\bf Q}_a$ is such that,
for given $\wt{\bf Q}_b$, the Frobenius norm of the first $R$ rows of
$\wt{\bf Q}_a^T\,[{\bf Z}_1{\bf Q}_b|\,{\bf Z}_2{\bf Q}_b]$ is maximal.

As in the proof of Theorem~\ref{t-8} (see (\ref{eq-eigenvec8})), when
$I>2R$ it also follows that ${\bf V}_1^T\wt{\bf Q}_a$ is of the form
(\ref{eq-eigenvecsuper359-1}).

Next, we consider $\wt{\bf Q}_b$ for given $\wt{\bf Q}_a$. We have
\be
\label{eq-unfold359-2}
\left[\begin{array}{c}
{\bf Q}_a^T{\bf Z}_1\\
{\bf Q}_a^T{\bf Z}_2\end{array}\right]\,\wt{\bf Q}_b=
\left[\begin{array}{cc}
{\bf R}_1 & \wt{\bf L}_1\\
{\bf R}_2 & \wt{\bf L}_2\end{array}\right]\,.
\ee

\noindent Let the SVD of $\left[\begin{array}{c}
{\bf Q}_a^T{\bf Z}_1\\
{\bf Q}_a^T{\bf Z}_2\end{array}\right]$ ($2R\times J$) be given by ${\bf
W}_2\,{\bf S}_2\,{\bf V}_2^T$, with ${\bf V}_2^T{\bf V}_2={\bf I}_J$ and
${\bf W}_2^T{\bf W}_2={\bf I}_{2R}$. Premultiplying (\ref{eq-unfold359-2})
by its transpose yields
\be
\label{eq-eigen359-2}
({\bf V}_2^T\wt{\bf Q}_b)^T\,{\bf S}_2^T{\bf S}_2\,({\bf V}_2^T\wt{\bf
Q}_b)= \left[\begin{array}{cc} {\bf K}_1 & {\bf O}_{R,J-R}\\ {\bf
O}_{J-R,R} & {\bf K}_2\end{array}\right]\,, \ee

\noindent where ${\bf K}_1$ is $R\times R$ and symmetric and
nonsingular, and ${\bf K}_2$ is $(J-R)\times (J-R)$ and symmetric. Matrix
${\bf S}_2$ is $2R\times J$ and contains the $\min(J,2R)$ nonzero singular
values on its diagonal.

Suppose $J\le 2R$. Then ${\bf S}_2^T{\bf S}_2$ is $J\times J$
diagonal and nonsingular. Since the diagonal entries of ${\bf S}_2^T{\bf
S}_2$ are positive and distinct, and ${\bf V}_2^T\wt{\bf Q}_b$ is orthonormal,
equation (\ref{eq-eigen359-2}) is the eigendecomposition of a $J\times J$
symmetric matrix that is nonsingular. In fact, (\ref{eq-eigen359-2}) can be
obtained as the superposition of the eigendecompositions of ${\bf K}_1$ and
${\bf K}_2$ (both nonsingular). This implies that
\be
\label{eq-eigenvecsuper359-2}
{\bf V}_2^T\wt{\bf Q}_b=\left[\begin{array}{cc}
\wt{\bf Q}_b^{(1)} & {\bf O}_{R,J-R}\\
{\bf O}_{J-R,R} & \wt{\bf Q}_b^{(2)}\end{array}\right]\,,
\ee

\noindent where $\wt{\bf Q}_b^{(1)}$ ($R\times R$) and $\wt{\bf Q}_b^{(2)}$
($(J-R)\times (J-R)$) are orthonormal. Hence, $\wt{\bf Q}_b$ is such that,
for given $\wt{\bf Q}_a$, the Frobenius norm of the first $R$ columns of
$\left[\begin{array}{c}
{\bf Q}_a^T{\bf Z}_1\\
{\bf Q}_a^T{\bf Z}_2\end{array}\right]\,\wt{\bf Q}_b$ is maximal. As above,
when $J>2R$ it also follows that ${\bf V}_2^T\wt{\bf Q}_b$ is of the
form (\ref{eq-eigenvecsuper359-2}).

From the above (also see the first-order optimality
conditions of problem (\ref{prob-W359}) in section 5.3), it follows that
$\wt{\bf Q}_a$ and $\wt{\bf Q}_b$ are such that the Frobenius norm of
$\wt{\bf G}_k={\bf R}_k$, $k=1,2$ is maximal in (\ref{eq-Zt359}).
Therefore, we have obtained an optimal solution ${\cal X}$ of problem
(\ref{prob-W359}) with slices ${\bf X}_k={\bf Q}_a\,{\bf R}_k\,{\bf
Q}_b^T$, $k=1,2$, and rank$({\cal X})=$ rank$({\cal R})=R$. This
contradicts the assumption in Theorem~\ref{t-359} that all optimal
solutions ${\cal X}$ of problem (\ref{prob-W359}) have rank larger than
$R$. Hence, $({\rm D}i)$ must hold.

Finally, we prove $({\rm D}ii)$ by contradiction. Suppose ${\bf R}_2{\bf
R}_1^{-1}$ has $R$ linearly independent eigenvectors. Hence, rank$({\cal
R})=R$ by Lemma~\ref{lem-rank} $(i)$. As in the proof of Theorem~\ref{t-1},
the optimality conditions (\ref{eq-stat12})--(\ref{eq-stat22}) then imply
that $\wt{\bf G}_k$, $k=1,2$, are upper triangular. As in the proof of
$({\rm D}i)$ above, we obtain an optimal solution ${\cal X}$ of problem
(\ref{prob-W359}) with rank$({\cal X})=$ rank$({\cal R})=R$, which is a
contradiction. Hence, $({\rm D}ii)$ must hold. This completes the proof of
Theorem~\ref{t-359}.

\end{document}